\documentclass[a4paper,12pt]{article}
\usepackage{amsmath,amsthm,verbatim,amssymb,amsfonts,amscd, graphicx}
\usepackage[all]{xy}
\usepackage{xcolor}
\usepackage{geometry}
\usepackage{graphicx}
\usepackage{hyperref}
\usepackage{enumitem}
\definecolor{darkblue}{rgb}{0,0,0.3}
\definecolor{urlblue}{rgb}{0,0,0.7}
\hypersetup{
	colorlinks=true,
	citecolor=blue,
	linkcolor=darkblue,
	urlcolor=urlblue,
}

\newcommand{\RR}{\mathbb{R}}

\DeclareMathOperator{\sgn}{sgn}
\DeclareMathOperator{\Area}{Area}
\DeclareMathOperator{\diam}{diam}

\newcommand{\D}{\nabla}
\newcommand{\p}{\partial}
\renewcommand{\th}{\theta}

\renewcommand{\bar}{\overline}
\renewcommand{\tilde}{\widetilde}
\renewcommand{\leq}{\leqslant}
\renewcommand{\geq}{\geqslant}
\renewcommand{\epsilon}{\varepsilon}

\newtheorem{theorem}{Theorem}[section]
\newtheorem{lemma}[theorem]{Lemma}
\newtheorem{cor}[theorem]{Corollary}
\newtheorem{prop}[theorem]{Proposition}
\newtheorem{defn}[theorem]{Definition}
\newtheorem*{fact}{Fact}
\numberwithin{equation}{section}
\theoremstyle{definition}

\theoremstyle{definition}
\newtheorem{remark}[theorem]{Remark}
\theoremstyle{definition}
\newtheorem{example}[theorem]{Example}

\DeclareMathOperator{\iso}{iso}
\DeclareMathOperator{\Ch}{Ch}
\DeclareMathOperator{\ID}{ID}
\DeclareMathOperator{\IN}{IN}

\begin{document}
	
\title{On closed surfaces with nonnegative curvature in the spectral sense}
\author{Kai Xu}
\date{}
\maketitle

\begin{abstract}
	We study closed orientable surfaces satisfying the spectral condition $\lambda_1(-\Delta+\beta K)\geq\lambda\geq0$, where $\beta$ is a positive constant and $K$ is the Gauss curvature. This condition naturally arises for stable minimal surfaces in 3-manifolds with positive scalar curvature. We show isoperimetric inequalities, area growth theorems and diameter bounds for such surfaces. The validity of these inequalities are subject to certain bounds for $\beta$. Associated to a positive super-solution $\Delta\varphi\leq\beta K\varphi$, the conformal metric $\varphi^{2/\beta}g$ has pointwise nonnegative curvature. Utilizing the geometry of the new metric, we prove H\"older precompactness and almost rigidity results concerning the main spectral condition.
\end{abstract}

\vspace{24pt}

\section{Introduction}

In this article, we study closed orientable surfaces that satisfy the spectral condition
\begin{equation}\label{eq-intro:eigenvalue_condition_main}
	\lambda_1(-\Delta+\beta K)\geq\lambda\geq0,
\end{equation}
where $K$ denotes the Gaussian curvature of $\Sigma$, and $\beta>0$ is a fixed constant. Some equivalent descriptions are
\begin{equation}\label{eq-intro:positivity_condition_main}
	\int_\Sigma\big(|\D\varphi|^2+\beta K\varphi^2\big)\,dA\geq\lambda\int_\Sigma\varphi^2\,dA,\qquad\forall\varphi\in C^\infty(\Sigma),
\end{equation}
and
\begin{equation}\label{eq-intro:exist_eigenfunction}
	\text{there exists }\varphi>0\text{ such that }\Delta\varphi\leq(\beta K-\lambda)\varphi.
\end{equation}
Clearly these conditions are implied by pointwise curvature bound $K\geq\beta^{-1}\lambda$. By setting $\varphi=1$ in (\ref{eq-intro:positivity_condition_main}), we obtain $2\pi\beta\chi(\Sigma)\geq\lambda|\Sigma|$, hence $\Sigma$ is topologically either a sphere or a torus. In fact, $\Sigma$ is either a sphere or a flat torus (Lemma \ref{lemma-cwt: perturbation}). We assume for the rest of this article that $\Sigma$ is topologically a sphere, for which interesting results are obtained. \\

One motivation for Condition (\ref{eq-intro:eigenvalue_condition_main}) is the study of positive scalar curvature (PSC) and general relativity. Let $\Sigma$ be a stable minimal surface in a 3-dimensional manifold $M$ with scalar curvature $R_M\geq R_0\geq0$. The second variation formula gives
\[\int_\Sigma\Big(|\D_\Sigma\varphi|^2+\frac12(2K_\Sigma-R_M-|h|^2)\varphi^2\Big)\,dA_\Sigma\geq0\quad (\forall\varphi),\]
and implies (\ref{eq-intro:positivity_condition_main}) with $\beta=1$ and $\lambda\geq R_0/2$.

The study of stable minimal surfaces is a key step for many important results in positive scalar curvature, such as Geroch's conjecture \cite{Schoen-Yau_1979a} and the positive mass theorem \cite{Schoen-Yau_1979c}. As black hole horizons can be mathematically described as minimal surfaces, condition (\ref{eq-intro:eigenvalue_condition_main}) is also related to topics in general relativity such as the Bartnik mass \cite{Li-Mantoulidis_2021, Mantoulidis-Schoen_2015}. For these background studies, we refer the reader to the book of Lee \cite{Lee}. Many new results in positive scalar curvature are recently proved using Gromov's $\mu$-bubbles \cite{Chodosh-Li, Gromov_2018, Gromov_2020_no_PSC, Lesourd-Unger-Yau_2021}, and we notice that a stable $\mu$-bubble (with suitable prescribed mean curvature functions) in a 3-manifold with uniformly positive scalar curvature also satisfies (\ref{eq-intro:eigenvalue_condition_main}) with $\beta=1$, $\lambda>0$. Extensive studies have been carried out for the non-compact case, see for example \cite{Berard-Castillon_2014, Castillon_2006, Colding-Minicozzi_2002, Fischer-Colbrie-Schoen_1980, Gromov-Lawson_1983, Munteanu-Sung-Wang_2023, Schoen_1983} and references therein. In higher dimensions there are two analogues of Condition (\ref{eq-intro:eigenvalue_condition_main}), in which one replaces  Gauss curvature either with scalar curvature \cite{Hirsch-Kazaras-Khuri-Zhang_2023, Li-Mantoulidis_2021}, or with the minimal eigenvalue of the Ricci curvature \cite{Carron_2019, Carron-Rose_2021, Xu_2023}. Analogues of Theorem \ref{thm-intro:isop_ineq}-\ref{thm-intro:bonnet_myers} can be obtained under the latter condition. In this paper we focus on the compact case in dimension two, for which the conformal structure is essential to many of our techniques.

\subsection{Geometric inequalities}

\vspace{9pt}

The initial motivation for this work is the question whether (\ref{eq-intro:eigenvalue_condition_main}) can be understood as a global positivity condition on curvature. Notice that (\ref{eq-intro:eigenvalue_condition_main}) gives no control on the local geometry: in \cite{Mantoulidis-Schoen_2015} Mantoulidis and Schoen constructed $C^1$-small but $C^2$-uncontrolled conformal perturbations of the round sphere, such that $\int \max(0,-K)\,dA$ is arbitrarily large. On the other hand, $\lambda_1(-\Delta+K)>0$ is preserved since it depends $C^1$-continuously on the conformal factor. The natural question then becomes asking what to expect globally. We give positive answers in terms of geometric inequalities. It turns out that (\ref{eq-intro:eigenvalue_condition_main}) implies isoperimetric inequality, volume comparison, and Bonnet-Myers' theorem (in the case $\lambda>0$).

\begin{defn}
	Given a closed Riemannian manifold $M$ of dimension $n$, the (Neumann) isoperimetric ratio is defined as
	\[\IN(M):=\inf_{\Omega\subset M}\frac{|\p\Omega|^n}{\min\big\{|\Omega|,|\Omega^c|\big\}^{n-1}}.\]
	The Cheeger constant is defined as
	\[\Ch(M):=\inf_{\Omega\subset M}\frac{|\p\Omega|}{\min\big\{|\Omega|,|\Omega^c|\big\}}.\]
	For a manifold with non-empty boundary, the (Dirichlet) isoperimetric ratio is defined as
	\[\ID(M):=\inf_{\Omega\subset\subset M}\frac{|\p\Omega|^n}{|\Omega|^{n-1}}.\]
\end{defn}

$ $

In analytic aspects, it is known that the isoperimetric ratio (resp. Cheeger constant) is related to the optimal constant in the Sobolev inequality (resp. Poincar\'e inequality).

The geometric inequalities we obtain are the following:

\begin{theorem}[isoperimetric inequality]\label{thm-intro:isop_ineq} {\ }
	
	Let $\Sigma$ be a closed surface that satisfies $\lambda_1(-\Delta+\beta K)\geq0$.
	
	(1) When $\beta>\frac12$, we have
	\begin{equation}\label{eq-intro:Ch_beta>1/2}
		\IN(\Sigma)\geq\frac{(2\beta-1)^2}{16\beta^2}\cdot\frac{|\Sigma|}{\diam(\Sigma)^2}.
	\end{equation}
	Moreover, the following local result holds: let $\gamma$ be a closed loop in $\Sigma$, and $N_\rho$ be its collar neighborhood of distance $\rho$. Assume $N_\rho$ is compact. If the first Dirichlet eigenvalue bound $\lambda_1^D(-\Delta+\beta K)\geq0$ is satisfied on $N_\rho$, then
	\[|N_\rho|\leq\frac{4\beta}{2\beta-1}\rho\,|\gamma|.\]
	
	(2) When $\frac12\geq\beta>\frac14$, for any $\epsilon>0$ we have
	\begin{equation}\label{eq-intro:IN_beta>1/4}
		\IN(\Sigma)\geq C(\beta,\epsilon)\big(\frac{|\Sigma|}{\diam(\Sigma)^2}\big)^{\frac1{4\beta-1}+\epsilon}.
	\end{equation}
\end{theorem}

It is not hard to see that $\Ch(\Sigma)^2\geq2\IN(\Sigma)/|\Sigma|$, so the Cheeger constant is bounded from below as well. In case (1) we have $\Ch(\Sigma)\geq\frac{2\beta-1}{2\sqrt 2\beta}\frac1{\diam(\Sigma)}$, which can be compared with Burago and Zalgaller's classical isoperimetric inequality $\Ch(\Sigma)\geq 1/\diam(\Sigma)$ under pointwise condition $K\geq 0$ \cite{Burago-Zalgaller}. When $\frac14<\beta\leq\frac12$, we cannot bound the Cheeger constant solely by the diameter, in contrast to the the case $\beta>\frac12$ and case of pointwise curvature bounds. See Example \ref{ex-countereg:collapsing} for details.

\begin{theorem}[volume comparison]\label{thm-intro:volume_comparison} {\ }
	
	Let $\beta>\frac14$ and $\Sigma$ be a closed surface with $\lambda_1(-\Delta+\beta K)\geq0$. Then there are universal constants $C=C(\beta,\frac{|\Sigma|}{\diam(\Sigma)^2})$ and $C'=C'(\beta)$, such that $Cr^2\leq|B(x,r)|\leq C'r^2$ for any geodesic ball of radius $r$. In particular, $|\Sigma|\leq C'(\beta)\diam(\Sigma)^2$.
\end{theorem}

The volume upper bound in Theorem \ref{thm-intro:volume_comparison} was essentially proved in Castillon \cite{Castillon_2006}, and the lower bound is a consequence of the isoperimetric inequality.

\begin{theorem}[weak Bonnet-Myers' Theorem]\label{thm-intro:bonnet_myers} {\ }
	
When $\beta>\frac14$ and $\lambda>0$, for a complete surface $\Sigma$ satisfying $\lambda_1(-\Delta+\beta K)\geq\lambda$ we have
\begin{equation}\label{eq-intro:bonnet-myers}
	\diam(\Sigma)\leq\frac{2\pi\beta}{\sqrt{\lambda(4\beta-1)}}.
\end{equation}
In particular, such surface must be compact.
\end{theorem}

The method for proving this statement is present in the literature, see \cite{Chodosh-Li, Gromov_2020_no_PSC, Schoen-Yau_1983, Shen-Ye_1996}. We give two proofs in Appendix A. Gromov-Lawson \cite{Gromov-Lawson_1983} proved a generalized statement for $\beta=1$ (the argument is easily generalized to $\beta>\frac12$), using the test function method to be introduced below. The author is informed that G. Carron and T. Richard has recently obtained the same diameter bound. Theorem \ref{thm-intro:bonnet_myers} is not true for $\beta\leq\frac14$; counterexamples include the hyperbolic disk for $\beta<\frac14$ (since $\lambda_1(-\Delta)=\frac14$), and the metric $g=dr^2+e^{4\lambda r^2}d\th^2$ with eigenfunction $\varphi=e^{-\lambda r^2}$ for $\beta=\frac14$.

\begin{remark}
	There are two thresholds $\beta=\frac12$ and $\beta=\frac14$ in the above theorems. The former one was observed by Gromov-Lawson \cite{Gromov-Lawson_1983} (see Proposition 8.11, Theorem 10.2). This is the critical value of $\beta$ for which local results can be obtained. By ``local'' we mean that (\ref{eq-intro:eigenvalue_condition_main}) with Dirichlet condition in a domain implies geometric inequalities in the same domain. When $\beta\in(\frac14,\frac12]$, the main geometric inequalities rely in an essential way on the closedness of $\Sigma$, and is no longer local in the above sense. The hyperbolic plane suggests that substantial change occurs to the geometry of (\ref{eq-intro:eigenvalue_condition_main}) below $\beta=\frac14$, and here we confirm the sharpness of $\beta=\frac14$ as a threshold. The same threshold also appeared in Castillon \cite{Castillon_2006}. It is possible to obtain further results for $\beta\leq\frac14$ under additional assumptions such as asymptotic area growth, see for example \cite{Espinar-Rosenberg_2011}.
\end{remark}

\subsection{The direct approach}

\vspace{9pt}

The main theorems are derived by testing (\ref{eq-intro:eigenvalue_condition_main}) with functions of the form $\varphi=\phi(d(x))$, where $d$ is the distance to a properly chosen smooth curve. Using coarea formula and integration by part, one transforms (\ref{eq-intro:eigenvalue_condition_main}) to an integral inequality involving $\phi(s)$ and $L(s):=|\{d(x,\gamma)=s\}|$, see (\ref{eq-direct:fundamental_eq}). We call this the fundamental equation. With different choice of $\phi$ applied to the fundamental equation, we obtain the main geometric inequalities. The core issue in this method is the low regularity of $L(s)$ due to the cut locus, and this was resolved through a careful analysis on the geodesic parallel coordinates \cite{Fiala_1940, Hartman_1964, Shiohama-Tanaka_1989, Shiohama-Tanaka_1993}. Still, this approach is restricted to dimension since the regularity issue is not resolved in higher dimensions. The fundamental equation is an efficient tool for studying condition (\ref{eq-direct:main_condition}). The earliest applications, to the author's knowledge, were due to Colding-Minicozzi \cite{Colding-Minicozzi_2002} and Gromov-Lawson \cite{Gromov-Lawson_1983}. For some other results obtained, we refer the reader to B\`erard-Castillon \cite{Berard-Castillon_2014}, Castillon \cite{Castillon_2006}, Munteanu-Sung-Wang \cite{Munteanu-Sung-Wang_2023}. See Section \ref{sec:direct} for details in the argument.

\subsection{Perspectives under a conformal change}

\vspace{9pt}

The main novelty of the present paper is the following conformal transformation technique, which provides us with new insights on (\ref{eq-intro:eigenvalue_condition_main}). Let $\varphi>0$ be any super-solution to $-\Delta+\beta K$, i.e. $\Delta\varphi\leq\beta K\varphi$. Define $u=\frac1\beta\log\varphi$, and we consider the conformal change $\tilde g=e^{2u}g$. The Gauss curvature of $\tilde g$ is computed to be
\begin{equation}\label{eq-intro:K>|du|^2}
	\tilde K\geq\beta|\tilde\D u|^2.
\end{equation}
Qualitively this is analogous to Schoen-Yau's observation that minimal surfaces in PSC manifolds are Yamabe positive. The new observation is the interplay between the geometry of $g$ and $\tilde g$. \\

The following anti-Harnack inequality a consequence of the isoperimetric inequality for $\tilde g$. In particular it applies to any counterexample of Theorem \ref{thm-intro:isop_ineq}(1) when $\beta\leq\frac12$. It is unlikely to obtain Harnack inequalities for $\varphi$, see Remark \ref{rmk-cwt:discussion}.
\begin{theorem}\label{thm-intro:anti_Harnack}
	Let $\beta>0$. Suppose $\Sigma$ is a closed surface, and $\varphi>0$ be any smooth function satisfying $\Delta\varphi\leq\beta K\varphi$. Then
	\[\frac{\max_\Sigma(\varphi)}{\min_\Sigma(\varphi)}\geq\big[\diam(\Sigma)\cdot\Ch(\Sigma)\big]^{-\beta/2}.\]
\end{theorem}

Since $(\Sigma,\tilde g)$ is a compact surface with nonnegative curvature, its geometry is controlled by the diameter-area ratio $\diam(\tilde g)^2/|\Sigma|_{\tilde g}$. We prove in Theorem \ref{thm-cwt:diam_bound} a uniform upper bound of this ratio, which is a result of the uniform Sobolev inequality for $g$. Based on this we show that
\begin{fact}
	$\tilde g$ is uniformly bi-Lipschitz equivalent to the round sphere after a suitable normalization.
\end{fact}
In this way $\tilde g$ can be understood as an approximate uniformization of $g$. The proof relies on Weyl's embedding theorem and therefore requires closedness of $\Sigma$. Finally, we obtain uniform H\"older equivalence with the standard metric for metrics satisfying (\ref{eq-intro:eigenvalue_condition_main}), as an implication of the theory of strong $A_\infty$ weights \cite{Aldana-Carron-Tapie_2021, David-Semmes_1990}. See Remark \ref{rmk-cwt:discussion} for further discussions in this direction.
\begin{theorem}\label{thm-intro:uniform-bi-Holder}
	Let $\beta>\frac14$. Given any constants $0<\alpha<1$ and $A_0>0$, there exists a constant $C=C(\alpha,\beta,A_0)$ such that: for any closed surface $\Sigma$ with $\diam(\Sigma)=1$, $|\Sigma|\geq A_0$ and satisfying (\ref{eq-intro:eigenvalue_condition_main}), there exists a diffeomorphism $\Psi:(\Sigma,g)\to(S^2,g_0)$, where $g_0$ is the metric with constant curvature $1$, such that $C^{-1}d(x,y)^{1/\alpha}\leq d_{g_0}(\Psi(x),\Psi(y))\leq Cd(x,y)^\alpha$.
\end{theorem}

Both Theorem \ref{thm-intro:volume_comparison} and \ref{thm-intro:uniform-bi-Holder} implies that the space of metrics with $\diam=1$, $\Area\geq A_0$ and satisfying (\ref{eq-intro:eigenvalue_condition_main}) is Gromov-Hausdorff precompact. \\

Finally, we consider convergence of metrics in the situation of almost rigidity. Observe that (\ref{eq-intro:eigenvalue_condition_main}) implies $|\Sigma|\leq 4\pi\beta\lambda^{-1}$, and equality case implies that $\Sigma$ has constant curvature $\lambda\beta^{-1}$ (since $\varphi=1$ must be the first eigenfunction). We show that $g$ metrically converges to the round sphere when $|\Sigma|$ approaches $4\pi\beta\lambda^{-1}$.
\begin{theorem}\label{thm-intro:almost_rigidity}
	Let $\beta>\frac14$. For any $\epsilon>0$ there is a constant $\delta=\delta(\epsilon,\beta)>0$ such that: if $\Sigma$ is topologically a sphere and satisfies $\lambda_1(-\Delta+\beta K)\geq\beta$, $|\Sigma|\geq 4\pi-\delta$, then there exists a diffeomorphism $\Psi:\Sigma\to(S^2,g_0)$ such that
	\[\max_{x,y\in\Sigma}\big|d_{g_0}(\Psi(x),\Psi(y))-d_g(x,y)\big|<\epsilon.\]
\end{theorem}

We have assumed $\lambda=\beta$ in the theorem; the general case only differs by a scaling. The proof consists of two steps. The first step is to show that $\tilde g$ (as defined above) bi-Lipschitz converges to the standard metric, and the second step is to show that $g$ converges to $\tilde g$ uniformly in the sense of distance function. See Subsection \ref{subsec:rig} for the detailed proof.

\vspace{12pt}

This paper is organized as follows. In Section \ref{sec:direct} we derive the fundamental equation (\ref{eq-direct:fundamental_eq}) and use it to prove Theorem \ref{thm-intro:isop_ineq}, \ref{thm-intro:volume_comparison}. In Section \ref{sec:aux} we prove several auxiliary lemmas to be used in the next section. In Section \ref{sec:cwt} we introduce the new perspective of conformal change mentioned above. The first subsection contains preparation works and the short proof of Theorem \ref{thm-intro:anti_Harnack}. The remaining two subsections \ref{subsec:bi-lipschitz}, \ref{subsec:rig} are devoted to the proof of Theorem \ref{thm-intro:uniform-bi-Holder} and \ref{thm-intro:almost_rigidity}. In Appendix \ref{sec:bonnet_myers} we give two proofs of Theorem \ref{thm-intro:bonnet_myers}. In Appendix \ref{sec:countereg} we discuss counterexamples for $\frac14<\beta\leq\frac12$. \\

\textbf{Notations.} We use $\Sigma$ to denote a closed surface, and $K$ to denote its Gaussian curvature. We use $|\cdot|$ to denote the area or length of an object, whose dimension is understood from the context. Volume forms are often suppressed when there is no ambiguity. We assume that all surfaces are smooth, connected, closed and orientable (for non-orientable surface, one can pass to its orientable double cover). We use $C$ to denote generic constants, which usually varies from term to term. The dependence of constants is indicated at the beginning of each section. \\

\textbf{Acknowledgements.} The author would like to thank his advisor Hubert Bray for encouragements during the progress of this work, as well as Sven Hirsch for helpful comments on the previous drafts of this paper. Also, the author would like to thank Simon Brendle, Gilles Carron, Demetre Kazaras, Marcus Khuri, Chao Li, Peter McGrath and Alec Payne for helpful conversations.

\section{The Method of Test Functions}\label{sec:direct}

\vspace{9pt}

Assume
\begin{equation}\label{eq-direct:main_condition}
	\int_\Sigma \big(|\D\varphi|^2+\beta K\varphi^2\big)\,dA\geq\lambda\int_\Sigma\varphi^2\,dA\quad(\forall\varphi\in C^\infty)
\end{equation}
on a closed surface $\Sigma$. In the first subsection we derive the fundamental equation (\ref{eq-direct:fundamental_eq}) from (\ref{eq-direct:main_condition}), extending the similar arguments in \cite{Castillon_2006, Gromov-Lawson_1983}. The equation is applicable to separating closed curves on closed surfaces, as suiting our needs for proving isoperimetric inequality. The condition $\phi''\geq0$ in \cite{Berard-Castillon_2014, Castillon_2006} is removed. We present a detailed derivation here for the sake of completeness. In the second subsection, we apply the fundamental equation to prove the main geometric inequalities. Linear functions are sufficient or the case $\beta>\frac12$, whereas functions with polynomial decay are used for the case $\frac14<\beta\leq\frac12$. Functions similar to the latter were used in \cite{Castillon_2006} to prove Euclidean volume growth for non-compact surfaces.

\subsection{The fundamental equation}\label{subsec:fundamental_eq}

\vspace{9pt}

Let $\Sigma$ be an oriented closed surface of Euler characterietic $\chi(\Sigma)$, and $\gamma\subset\Sigma$ be a connected smooth closed curve. For the moment we assume that $\gamma$ is separating, so that $\Sigma\setminus\gamma=\Omega^+\cup\Omega^-$. Define the signed distance function
\[d(x):=\left\{\begin{aligned}
	& d(x,\gamma)\qquad(x\in\Omega^+), \\
	& -d(x,\gamma)\qquad(x\in\Omega^-).
\end{aligned}\right.\]
Let $-\rho^-:=\min(d)$, $\rho^+:=\max(d)$. Define
\begin{equation}\label{eq-direct:definitions1}
	\gamma(s):=\{d=s\},\qquad\Omega(s):=\left\{\begin{aligned}
		& \{0\leq d< s\}\quad(s\geq0), \\
		& \{s>d\leq 0\}\quad(s\leq 0),
	\end{aligned}\right.\qquad\chi(s):=\chi(\Omega(s)),
\end{equation}

\begin{equation}\label{eq-direct:definitions2}
	L(s):=|\gamma(s)|,\qquad
	K(s):=\int_{\gamma(s)}K\,dl,\qquad
	G(s):=\int_0^s K(t)\,dt.\qquad
\end{equation}
Finally, define
\[\Gamma(s):=2\pi\chi(s)\sgn(s)-G(s)+\int_{\gamma}\kappa,\]
where $\kappa$ is the geodesic curvature of $\gamma$. We adopt the sign convention that the unit normal vector always points into $\Omega^+$. Whenever $\gamma(s)$ is smooth, it follows from Gauss-Bonnet formula that $\Gamma(s)=\int_{\gamma(s)}\kappa$. Also, $L'(s)=\Gamma(s)$ when $|s|$ is small. However, $L(s)$ in general not continuous because of cut locus (see Figure 1 in \cite{Hartman_1964}). The following partial regularity result can be obtained:
\begin{lemma}[\cite{Hartman_1964, Shiohama-Tanaka_1989, Shiohama-Tanaka_1993}]\label{lemma-direct:L'<Gamma}
	For any smooth closed curve $\gamma\subset M$, $\gamma(s)$ as defined above is piecewise smooth for almost every $s$. Moreover, $L(s)$ is almost everywhere differentiable and we have $L'(s)\leq\Gamma(s)$ for $s>0$, $L'(s)\geq\Gamma(s)$ for $s<0$. Finally, $L(s_2)\leq L(s_1)+\int_{s_1}^{s_2}L'(t)\,dt$ for a.e. $0\leq s_1<s_2\leq\rho^+$ and the inequality is reversed when $s_1<s_2<0$.
\end{lemma}

\begin{remark}
	In the case where $g$ and $\gamma$ are both real analytic, Fiala \cite{Fiala_1940} proved that $L(s)$ is continuous and real analytic except at finitely many times.
\end{remark}

This lemma implies that $L'\leq\Gamma$ distributionally when $s>0$ (and $L'\geq\Gamma$ when $s<0$). Now let $\phi:[-\rho^-,\rho^+]\to[0,\infty]$ be a non-negative continuous function, smooth except at $s=0$, with $\phi'\leq0\,(\forall s>0)$ and $\phi'\geq0\,(\forall s<0)$. We test (\ref{eq-direct:main_condition}) with function $\phi(d(x))$, by coarea formula we have:
\begin{equation}\label{eq-direct-1}
	\lambda\int_{-\rho^-}^{\rho^+}L(s)\phi(s)^2\,ds\leq\int_{-\rho^-}^{\rho^+}\Big[L(s)\phi'(s)^2+\beta K(s)\phi(s)^2\Big]\,ds.
\end{equation}
From Gauss-Bonnet formula and integration by part, we obtain
\[\begin{aligned}
	\int_{-\rho^-}^{\rho^+}K\phi^2\,ds &=
	\big[G\phi^2\big]_{\rho^-}^{\rho^+} - \int_{-\rho^-}^{\rho^+} 2G\phi\phi' \\
	&= 2\pi\Big[\chi(\Omega^+)\phi(\rho^+)^2 + \chi(\Omega^-)\phi(\rho^-)^2\Big]
	+4\pi\Big(\int_{-\rho^-}^0\chi\phi\phi' - \int_0^{\rho^+}\chi\phi\phi'\Big) \\
	&\qquad +\int_{-\rho^-}^{\rho^+}2\Gamma\phi\phi'\,ds.
\end{aligned}\]

Since $\Omega(s)$ is connected with more than one boundary components, we have $\Omega(s)\leq0$ hence the integrals on $\chi\phi\phi'$ is non-positive. For the $\int \Gamma\phi\phi'$ term, we claim that
\[\int_0^{\rho^+}\Gamma f\,ds\leq-\int_0^{\rho^+}Lf'\,ds-L(0)f(0)\]
for any non-positive $f\in C^\infty([0,\rho^+])$. Since $L'=\Gamma$ is smooth in a neighborhood of zero, it suffices to prove the identity for $f$ supported away from zero. Choose $\epsilon\ll\delta\ll1$ such that $L$ is differentiable at $\rho^+-\delta$. Choose a cutoff function $\psi$ with $\psi|_{[0,\rho^+-\delta-\epsilon]}\equiv0$, $\psi|_{[\rho^+-\delta+\epsilon,\rho^+]}\equiv1$, and $|\psi'|\leq 4\epsilon^{-1}$. Decompose $f=f\psi+f(1-\psi)$. Note that $\int\Gamma f\psi\,ds\to0$ as $\delta\to0$. Since $L'\leq\Gamma$ weakly, we have
\[\int_0^{\rho^+}\Gamma f(1-\psi)\,ds\leq-\int_0^{\rho^+} L\big(f(1-\psi)\big)'\,ds=-\int_0^{\rho^+} Lf'\,ds+o_\delta(1)+\int_0^{\rho^+} Lf\psi'\,ds.\]
Finally, $\int Lf\psi'\,ds\to L(\rho^+-\delta)f(\rho^+-\delta)$ as $\epsilon\to0$, and the latter is non-positive. This proves our claim.

For the same reason, we have $\int_{-\rho^-}^0\Gamma f\,ds\leq -\int_0^{\rho^+}Lf'\,ds+L(0)f(0)$ for all $f\geq0$. Applying $f=2\phi\phi'$ and combining all the inequalities, we obtain

\begin{equation}\label{eq-direct:fundamental_eq}
	\begin{aligned}
		&(2\beta-1)\int_{-\rho^-}^{\rho^+}L(\phi')^2\,ds
			+2\beta\int_{-\rho^-}^{\rho^+}L\phi\phi''\,ds
			+\lambda\int_{-\rho^-}^{\rho^+}L\phi^2\,ds \\
		&\qquad\qquad\leq
			2\pi\beta\Big[\chi(\Omega^+)\phi(\rho^+)^2+\chi(\Omega^-)\phi(\rho^-)^2\Big]
			+2\beta|\gamma|\phi(0)\big[\phi'_-(0)-\phi'_+(0)\big].
	\end{aligned}
\end{equation}

We summarize all the conditions into the following lemma:

\begin{lemma}\label{lemma-direct:separating}
	Assume that $\Sigma$ satisfies (\ref{eq-direct:main_condition}), and $\gamma\subset\Sigma$ is a connected separating smooth curve. Suppose $\phi:[-\rho^-,\rho^+]\to[0,\infty]$ is nonnegative and continuous, with $\phi'\leq0$ when $s>0$ and $\phi'\geq0$ when $s<0$. Then (\ref{eq-direct:fundamental_eq}) holds. For notations see (\ref{eq-direct:definitions1}) (\ref{eq-direct:definitions2}). \hfill$\Box$
\end{lemma}

\subsection{Derivation of geometric inequalities}

\vspace{9pt}

\begin{proof}[Proof of Theorem \ref{thm-intro:isop_ineq}] {\ }

Suppose $\Sigma=\Omega^+\cup\Omega^-$, with common smooth boundary $\gamma=\p\Omega^+=\p\Omega^-$. It suffices to prove the theorem when $\Omega^-$ is connected. To see this, suppose $\Omega^-=\Omega_1\cup\Omega_2$ as a disjoint union. If both $|\Omega_1|$ and $|\Omega_2|$ are less than $\frac12|\Sigma|$, then by induction on the number of connected components, we have $|\p\Omega_i|\geq|\Omega_i|^s C_{\iso}^s$ where $C_{\iso}=\IN(\Sigma)$ or $\Ch(\Sigma)$, and $s=1/2$ or $1$. Hence $|\p\Omega|\geq|\Omega|^sC_{\iso}^s$, which proves the isoperimetric inequality for $\Omega$. If $|\Omega_1|\geq\frac12|\Sigma|$, we have $|\p\Omega|\geq|\p\Omega_1|\geq C_{\iso}^s|\Sigma\setminus\Omega_1|^s\geq C_{\iso}^s|\Sigma\setminus\Omega|^s$.

Therefore, we assume that both $\Omega^+$ and $\Omega^-$ are connected, thus $\gamma$ is connected (since $\pi_1(\Sigma)=0$). To prove (1), we test equation (\ref{eq-direct:fundamental_eq}) with
\[\phi(s)=\left\{\begin{aligned}
	& (\rho^-+s)/\rho^-\qquad(s\leq0), \\
	& (\rho^+-s)/\rho^+\qquad(s\geq0).
\end{aligned}\right.\]
We obtain
\[(2\beta-1)\Big(\frac{|\Omega^-|}{(\rho^-)^2}+\frac{|\Omega^+|}{(\rho^+)^2}\Big)\leq 2\beta\,|\gamma|\,\Big(\frac{1}{\rho^-}+\frac{1}{\rho^+}\Big)\]

We may assume $\rho^-\leq\rho^+$. The first part of Theorem \ref{thm-intro:isop_ineq}(1) follows from
\[(2\beta-1)^2\frac{|\Omega^-|}{(\rho^-)^2}\cdot \frac{|\Sigma|}{\diam(\Sigma)^2}
	\leq 16\beta^2\frac{|\gamma|^2}{(\rho^-)^2}.\]
For the second part of (1), we choose $\rho^+=\rho^-=\rho$, and note that $\phi$ is indeed a test function for the Dirichlet eigenvalue condition (i.e. $\phi$ vanishes on $\p N_\rho$).

\vspace{18pt}

To prove (2), we test with the following function:
\begin{equation}\label{eq-direct:phi_polynomial}
	\phi(s)=\left\{\begin{aligned}
		& \frac{\big[(1+\sigma)\rho^-+s\big]^p}{(1+\sigma)^p(\rho^-)^p}\qquad(s\leq0), \\
		& \frac{\big[(1+\sigma)\rho^+-s\big]^p}{(1+\sigma)^p(\rho^+)^p}\qquad(s\geq0), \\
	\end{aligned}\right.
\end{equation}
where the coefficients $\sigma>0,p>1$ are to be chosen later. Equation (\ref{eq-direct:fundamental_eq}) now gives
\begin{equation}\label{eq-direct-3}
	I^-+I^+\leq 4\pi\beta\frac{\sigma^{2p}}{(1+\sigma)^{2p}}+\frac{2p\beta}{1+\sigma}\,|\gamma|\,\big(\frac1{\rho^-}+\frac1{\rho^+}\big),
\end{equation}
where
\begin{equation}\label{eq-direct-4}
	I^-:=\big[(4\beta-1)p^2-2\beta p\big]\frac1{(1+\sigma)^{2p}(\rho^-)^{2p}}\int_{-\rho^-}^0 L\big[(1+\sigma)\rho^-+s\big]^{2p-2}\,ds,
\end{equation}
and $I^+$ is defined analogously. The coefficient $(4\beta-1)p^2-2\beta p$ must be positive, thus $p>\frac{2\beta}{4\beta-1}$. We choose $p=\frac{2\beta+\epsilon}{4\beta-1}$, $\epsilon>0$. Now (\ref{eq-direct-3}) (\ref{eq-direct-4}) implies
\begin{equation}\label{eq-direct:before_chossing_sigma}
	C_1\frac{\sigma^{2p-2}}{(1+\sigma)^{2p}}\Big[\frac{|\Omega^-|}{(\rho^-)^2}+\frac{|\Omega^+|}{(\rho^+)^2}\Big]\leq 4\pi\beta\frac{\sigma^{2p}}{(1+\sigma)^{2p}}+\frac{2p\beta}{1+\sigma}\,|\gamma|\,\big(\frac1{\rho^-}+\frac1{\rho^+}\big),
\end{equation}
where $C_1=C_1(\beta,\epsilon)$ is a constant depending only on $\beta,\epsilon$. For convenience, we denote $Z:=\frac{|\Omega^-|}{(\rho^-)^2}+\frac{|\Omega^+|}{(\rho^+)^2}$. Choose $\sigma=\sqrt{\frac{C_1}{8\pi\beta}Z}$, for this choice we have
\begin{equation}\label{eq-direct:plugging_in_sigma}
	\frac{2p\beta}{1+\sigma}\,|\gamma|\,\big(\frac1{\rho^-}+\frac1{\rho^+}\big)\geq\frac{C_1}2\frac{\sigma^{2p-2}}{(1+\sigma)^{2p}}Z.
\end{equation}
Without loss of generality, assume $\rho^-\leq\rho^+$. Then we have
\[\frac{|\gamma|}{\rho^-}\geq C(\beta,\epsilon)\frac{\sigma^{2p-1}}{(1+\sigma)^{2p-1}}\sqrt Z\geq C(\beta,\epsilon)\min\big(1,\sigma^{2p-1}\big)\frac{\sqrt{|\Omega^-|}}{\rho^-},\]
or equivalently,
\begin{equation}\label{eq-direct:weaker_isop_ineq}
	\frac{|\gamma|^2}{|\Omega^-|}\geq C(\beta,\epsilon)\min\big(1,Z^{2p-1}\big)\geq C(\beta,\epsilon)\min\big(1,(\frac{|\Sigma|}{\diam(\Sigma)^2})^{\frac{1+2\epsilon}{4\beta-1}}\big).
\end{equation}
Combined with total volume upper bound to be proved below, we obtain the strongest form (\ref{eq-intro:IN_beta>1/4}).
\end{proof}

\begin{proof}[Proof of Theorem \ref{thm-intro:volume_comparison}] {\ }
	
	The volume upper bound was essentially proved in \cite[Proposition 2.2]{Castillon_2006}. Here we include a proof for the reader's convenience. Let $\epsilon\ll1$ such that the geodesic sphere $\gamma=\p B(x,\epsilon)$ is smooth. Consider the test function
	\[\phi(s)=\left\{\begin{aligned}
		& 1 \qquad (s\leq0) \\
		& (1-\frac{s}{2r})^p \qquad (0\leq s\leq 2r) \\
		& 0 \qquad (s\geq2r)
	\end{aligned}\right.\]
	Applying (\ref{eq-direct:fundamental_eq}) to $\gamma$, similar calculation yields
	\[\frac{C_1}{(2r)^{2p}}\int_0^{2r}L(2r-s)^{2p-2}\,ds\leq 4\pi\beta+2p\beta|\gamma|r^{-1}.\]
	The left hand side is greater than $Cr^{-2}|B(x,r+\epsilon)|$. Letting $\epsilon\to0$ we obtain $|B(x,r)|\leq Cr^2$. In particular, $|\Sigma|\leq C\diam(\Sigma)^2$. This proves the upper bound and shows that it depends only on $\beta$.
	
	The volume lower bound is a consequence of the isoperimetric inequality. let $r_1=\sup\{r: B(x,r)\leq\frac12|\Sigma|\}$. For any $r\leq r_1$ such that $\p B(x,r)$ is piecewise smooth, Theorem \ref{thm-intro:isop_ineq} gives
	\[|\p B(x,r)|\geq C(\beta,\frac{|\Sigma|}{\diam(\Sigma)^2})\sqrt{|B(x,r)|}.\]
	By coarea formula and gronwall's inequality we obtain $|B(x,r)|\geq C(\beta,\frac{|\Sigma|}{\diam(\Sigma)^2})r^2$. For $r\geq r_1$ we have $|\p B(x,r)|\geq\frac12|\Sigma|\geq\frac12\frac{|\Sigma|}{\diam(\Sigma)^2}\cdot r^2$. This completes the proof.
\end{proof}

\begin{remark}
	 In the non-compact case, the proof shows that for any $\delta>0$, $|B(x,r)|/r^2$ is uniformly bounded from above if $\lambda_1^D(-\Delta+\beta K)\geq0$ in $B(x,(1+\delta)r)$. When $\beta>\frac12$, the linear function $\phi(s)=1-s/r$ can be used to show that $|B(x,r)|\leq Cr^2$ whenever the metric ball $B(x,r)$ satisfies $\lambda_1^D(-\Delta+\beta K)\geq0$.
\end{remark}

\section{Auxiliary Lemmas}\label{sec:aux}

\vspace{9pt}

\subsection{Collection of facts in convex geometry}

\vspace{9pt}

\begin{lemma}\label{lemma-aux: convex geometry}
	
	(1) (Weyl's embedding theorem) Given any metric $(S^2,g)$ with curvature $K>0$, there exists an isometric embedding into $\RR^3$, unique up to rigid motion. The image of the embedding is the boundary of a strictly convex set.
	
	(2) Suppose $\Omega_1\subset\Omega_2$ are two domains in $\RR^n$, with $\Omega_1$ smooth convex and $\Omega_2$ piecewise smooth. Then the orthogonal projection from $\p\Omega_2$ to $\p\Omega_1$ is 1-Lipschitz.
	
	(3) For a closed surface $\Sigma$ with curvature $K\geq0$, we have $\Ch(\Sigma)\geq1/\diam(\Sigma)$ and $\IN(\Sigma)\geq |\Sigma|/\diam(\Sigma)^2$.
	
\end{lemma}

\begin{proof}
	(1) is a classical theorem, see Nirenberg \cite{Nirenberg_1953} and Pogorelov \cite{Pogorelov}.
	
	(2) is proved by computing the differential of the projection map. Let $\Phi:\p\Omega_2\to\p\Omega_1$ be the (well-defined) orthogonal projection map. Denote by $N, A$ the outer unit normal vector and second fundamental form of $\p\Omega_1$. Let $h(x)=d(x,\Phi(x))\geq0$. For any $x\in\p\Omega_2$ we have $x-\Phi(x)=h(x)N(\Phi(x))$. Differentiating this identity at a tangent vector $v$, we obtain $v-d\Phi_x(v)=dh_x(v)N(\Phi(x))+h(x)dN_{\Phi(x)}(d\Phi_x(v))$. Taking inner product with $d\Phi_x(v)$ we obtain $|v|\cdot|d\Phi_x(v)|-|d\Phi_x(v)|^2\geq\big(v-d\Phi_x(v)\big)\cdot d\Phi_x(v)=h(x)A(d\Phi_x(v),d\Phi_x(v))\geq0$. This shows $d\Phi$ is non-expanding.
	
	(3) Let $\Omega$ be any domain. An inequality of Burago-Zalgaller (see \cite{Burago-Zalgaller, Osserman_1978}) states that $\rho|\p\Omega|\geq|\Omega|+(\pi-\frac12\int_\Omega K)\rho^2$, where $\rho$ is the radius of the largest metric ball contained in $\Omega$. By possibly switching between $\Omega$ and $\Omega^c$, we may assume $\pi-\frac12\int_\Omega K\geq0$, hence $|\Omega|\leq\rho|\p\Omega|$. This gives the Cheeger constant bound. By relative volume comparison we have $|\p\Omega|^2\geq|\Omega|\cdot\frac{|\Omega|}{\rho^2}\geq|\Omega|\cdot\frac{|\Sigma|}{\diam(\Sigma)^2}$, which gives the isoperimetric ratio bound.
\end{proof}

\subsection{A precise form of Moser-Trudinger's inequality}

\vspace{9pt}

The classical Moser-Trudinger Sobolev inequality \cite{Moser_1971} states that for a domain $\Omega\subset\RR^2$ and a function $u\in H^1_0(\Omega)$ we have
\begin{equation}\label{eq-aux-1}
	\int_{\Omega}\exp\big(\frac{4\pi u^2}{\int_{\Omega}|\D u|^2}\big)\leq C|\Omega|
\end{equation}
for a universal constant $C$. For general domains in smooth surfaces, (\ref{eq-aux-1}) continue to hold with the same critical exponent $4\pi$, while the constant $C$ depends on the domain under consideration. When the surface has a conic singularity of cone angle $\th<2\pi$, the critical exponent becomes $2\th$ \cite{Chang-Yang_1988, Troyanov_1991}. We concern the case where only a lower bound on the isoperimetric ratio is known. In this case the round cone represents the worst control, as seen by P\'olya-Szeg\"o symmetrization. The following lemma is useful for later applications.

\begin{theorem}\label{thm-aux:M-T-Dirichlet}
	Let $(\Omega,g)$ be a smooth domain with non-empty boundary, such that $\ID(\Omega)\geq\xi$. Then for any function $u\in H^1_0(\Omega)$ we have
	\begin{equation}\label{eq-aux-2}
		\int_\Omega \exp\big(\frac{\xi u^2}{\int_\Omega|\D u|^2}\big)\leq C(\xi)\cdot|\Omega|,
	\end{equation}
	where $C(\xi)$ means a constant depending on $\xi$.
\end{theorem}
\begin{proof}
	We can replace $u$ by $|u|$ and therefore assume $u\geq0$. Since (\ref{eq-aux-2}) is scale-invariant, we may assume $|\Omega|=1$. We apply a P\'olya-Szeg\"o symmetrization procedure, comparing $\Omega$ to a model space that has smaller isoperimetric ratio. Let $(r,\th)$ be the polar coordinates on a disk $D^2$. Equip $D^2$ with the cone metric $g_0=dr^2+\epsilon^2 r^2d\th^2$ ($0\leq r\leq L, 0\leq\th\leq 2\pi$), where $\epsilon$ is determined by requiring $|\p\Omega|^2=\xi|\Omega|$ for concentric cones $\Omega$, and $L$ is determined such that $|D^2|_{g_0}=1$. Consider the radial function $u_0(r)\in H^1_0(D^2)$ uniquely determined by $du_0/dr\leq0$ and the condition
	\[\big|\{y\in D^2: u_0(y)>t\}\big|_{g_0}=\big|\{x\in\Sigma: u(x)>t\}\big|_g,\quad\forall t\geq0.\]
	From isoperimetric ratio comparison we have
	\[\big|\{y\in D^2: u_0(y)=t\}\big|_{g_0}\leq\big|\{x\in\Sigma: u(x)=t\}\big|_g.\]
	The standard argument (see for example \cite[Section III.1]{Schoen-Yau_lectures}), using rearrangement and the coarea formula, gives $\int_{D^2}\exp(pu_0^2)\,dA_0=\int_\Omega\exp(pu^2)\,dA$ and $\int_{D^2}|\D_{g_0}u_0|^2\,dA_0\leq\int_\Omega|\D u|^2\,dA$. The desired result follows from the Moser-Trudinger inequality with conic singularity on $(D^2,g_0)$ and $u_0$.
\end{proof}

\begin{cor}\label{cor-aux: M-T inequality for e^2u}
	For the same conditions as in Theorem \ref{thm-aux:M-T-Dirichlet}, we have
	\[\int_\Omega e^{2u}\leq C(\xi)\,|\Omega|\exp\big(\frac1\xi\int_\Omega|\D u|^2\big).\]
\end{cor}
\begin{proof}
	Combine Theorem \ref{thm-aux:M-T-Dirichlet} with Young's inequality $2u\leq \frac{\xi u^2}{\int|\D u|^2}+\frac1\xi\int|\D u|^2$.
\end{proof}
Corollary \ref{cor-aux: M-T inequality for e^2u} also applies to functions $u$ that are non-positive on the boundary, by replacing $u$ with $\max(u,0)$. Similar inequality also holds for closed manifold with Neumann isoperimetric bounds.

\begin{theorem}\label{thm-aux:M-T-closed}
	Let $(\Sigma,g)$ be a closed surface with $\IN(\Sigma)\geq\xi>0$. Then for any smooth function $u$ with $\int_\Sigma u=0$, we have
	\begin{equation}\label{eq-aux:IN-M-T-e^u^2}
		\int_\Sigma\exp\big(\frac{\xi u^2}{\int_\Sigma|\D u|^2}\big)\leq C(\xi)\cdot|\Sigma|.
	\end{equation}
	Hence
	\begin{equation}\label{eq-aux:IN-M-T-e^pu}
		\int_\Omega e^{pu}\leq C(\xi,p)\,|\Omega|\,\exp\big(\frac1\xi\int_\Sigma|\D u|^2\big).
	\end{equation}
\end{theorem}
\begin{proof}
	After a perturbation, we may assume that there exists a regular value $b$ such that $\big|\{u\geq b\}\big|=\big|\{u\leq b\}\big|=\frac12|\Sigma|$. Assume without loss of generality that $b\geq0$. By Chebyshev and Cheeger's inequality,
	\[\begin{aligned}
		\frac12|\Sigma| &= \big|\{u\geq b\}\big|\leq b^{-2}\int_\Sigma u^2\leq 4b^{-2}\Ch(\Sigma)^{-2}\int_\Sigma|\D u|^2 \\
		&\leq 2b^{-2}|\Sigma|\cdot\IN(\Sigma)^{-1}\int_\Sigma|\D u|^2.
	\end{aligned}\]
	Hence $b^2\leq 4\xi^{-1}\int|\D u|^2$. Apply Theorem \ref{thm-aux:M-T-Dirichlet} to $\Omega=\{u\geq b\}$ (note that $\ID(\Omega)\geq\IN(\Sigma)$):
	\[\int_\Omega\exp\big(\frac{\xi u^2}{\int|\D u|^2}\big)\leq \int_\Omega\exp\big(\frac{\xi (u-b)^2+\xi b^2}{\int|\D u|^2}\big)\leq C(\xi)\,|\Omega|.\]
	The same conclusion holds for $\{u\leq b\}$, therefore (\ref{eq-aux:IN-M-T-e^u^2}) holds. (\ref{eq-aux:IN-M-T-e^pu}) follows from (\ref{eq-aux:IN-M-T-e^u^2}) and Young's inequality.
\end{proof}

\section{The Perspective under Conformal Change}\label{sec:cwt}

\vspace{9pt}

\subsection{Preparation and first consequences}\label{subsec:preparation}

\vspace{9pt}

In this section we assume that $\Sigma$ is a closed surface satisfying (\ref{eq-intro:eigenvalue_condition_main}) and is topologically a sphere. Let $\varphi=e^{\beta u}$ satisfy $\Delta\varphi\leq\beta K\varphi$, with $\beta>\frac14$. We have
\begin{equation}\label{eq-cwt:supersolution}
	e^{\beta u}(\beta\Delta u+\beta^2|\D u|^2)=\Delta e^{\beta u}\leq\beta Ke^{\beta u}.
\end{equation}
Consider the conformal change $\tilde g=e^{2u}g$. We may normalize $\varphi$ so that $|\Sigma|_{\tilde g}=\int_\Sigma e^{2u}\,dA=1$. The Gauss curvature of $\tilde g$ is computed to be
\begin{equation}\label{eq-cwt:K>e^-2u+du^2}
	\tilde K=e^{-2u}(K-\Delta u)\geq e^{-2u}\beta|\D u|^2=\beta|\tilde\D u|^2,
\end{equation}
where we denote $|\tilde\D u|^2=|\tilde\D u|_{\tilde g}^2$ for brevity. The two basic observations from (\ref{eq-cwt:K>e^-2u+du^2}) are $\tilde K\geq0$ and $\int|\tilde\D u|^2\,d\tilde A\leq 4\pi\beta^{-1}$. The former gives geometric control on $\tilde g$, while the latter provides strong relation between $g$ and $\tilde g$. \\

\begin{proof}[Proof of Theorem \ref{thm-intro:anti_Harnack}] {\ }
	
	Suppose $\max_\Sigma(\varphi)=A\min_\Sigma(\varphi)$. Applying Burago-Zalgaller's isoperimetric inequality to $(\Sigma,\tilde g)$, we have
	\[\inf_{\Omega\subset\Sigma}\frac{\int_{\p\Omega}\varphi^{1/\beta}\,dl}{\min\big\{
		\int_\Omega\varphi^{2/\beta}\,dA,\int_{\Omega^c}\varphi^{2/\beta}\,dA\big\}}\geq\frac1{\diam(\Sigma,\tilde g)}.\]
	
	Note that the left hand side is $\leq\Ch(\Sigma)A^{2/\beta}(\max_\Sigma\varphi)^{-1/\beta}$, while the right hand side is $\geq(\max_\Sigma\varphi)^{-1/\beta}/\diam(\Sigma)$. Hence
	\[A^{2/\beta}\geq\frac{1}{\Ch(\Sigma)\cdot\diam(\Sigma)}.\]
	This proves the theorem.
\end{proof}

\begin{remark}\label{eq-cwt:no_Harnack}
	We do not expect global semi-Harnack inequalities for $\varphi>0$ satisfying $\Delta\varphi\leq\beta K\varphi$, for the following reason. First we observe that no uniform bound on $\sup(\varphi)$ can be found. Let $\varphi>0$ uniquely solve $\Delta\varphi=\beta\varphi-\delta_x$ on the round sphere. Applying the heat operator $P_\epsilon$, we obtain a family of smooth functions $\varphi_\epsilon$ with $\Delta\varphi_\epsilon\leq\beta\varphi_\epsilon$. Next, we note the following interesting relation:
	\[\Delta_g\varphi\leq\beta K_g\varphi\ \Rightarrow\ 
		\Delta_{g'}\varphi^{-1}\leq\beta K_{g'}\varphi^{-1},
		\quad\text{where }g'=\varphi^{4/\beta}g.\]
	This can be seen by reversing the implication from (\ref{eq-cwt:supersolution}) to (\ref{eq-cwt:K>e^-2u+du^2}) with $u\to-u$. Thus we have $\Delta_{g_\epsilon}\varphi_{\epsilon}^{-1}\leq\beta K_{g_\epsilon}\varphi_\epsilon^{-1}$ where $g_\epsilon=\varphi_\epsilon^{4/\beta}g$. Note that $\inf(\varphi_\epsilon^{-1})\to0$ as $\epsilon\to0$, whereas the diameter and area for $g_\epsilon$ can be verified to have uniform controls. This shows that uniform lower bounds on $\inf(\varphi)$ are unlikely to be obtained.
\end{remark}

The rest of this section relies on Weyl's embedding theorem for $(\Sigma,\tilde g)$, hence requires $\tilde K>0$. However, (\ref{eq-cwt:K>e^-2u+du^2}) only gives $\tilde K\geq0$. With the following lemma available, we may assume after a $C^2$ approximation that (\ref{eq-cwt:supersolution}) is a strict inequality, hence $\tilde K>\beta|\tilde\D u|^2$.

\begin{lemma}\label{lemma-cwt: perturbation}
	Let $(\Sigma,g)$ be a closed surface with $\lambda_1(-\Delta+\beta K)=0$, and let $\varphi$ be the first eigenfunction. If $\Sigma$ is not a flat torus, then for any $\epsilon>0$ there exists another metric $g'$ such that $||g'-g||_{C^2(g)}<\epsilon$ with $\lambda_1(-\Delta_{g'}+\beta K')>0$. Furthermore, the first eigenfunction $\varphi'$ satisfies $||\varphi'-\varphi||_{C^2(g)}<\epsilon$.
\end{lemma}
\begin{proof}
	Consider a smooth family of metrics $g(t)=e^{2f(t)}\bar g$, $f(0)=0$, with initial variation $\frac{df}{dt}|_{t=0}=h$. By \cite[Lemma A.1]{Mantoulidis-Schoen_2015}, the first eigenfunctions $\varphi_t$ of $-\Delta_{g(t)}+\beta K_{g(t)}$ constitute a smooth family. We may normalize so that $||\varphi_t||_{L^2(g_t)}=1$. We compute
	\begin{align}
		\frac{d\lambda_1}{dt}\Big|_{t=0} &=\frac{d}{dt}\Big|_{t=0}\int_{\Sigma}\Big(|\D_t\varphi_t|^2+\beta K_t\varphi_t^2\Big)\,dA_t \nonumber \\
		&= \frac{d}{dt}\Big|_{t=0}\int_{\Sigma}\Big(|\D_t\varphi|^2+\beta K_t\varphi^2\Big)\,dA_t \nonumber \\
		&= -\beta\int_\Sigma\Delta h\cdot\varphi^2\,dA \label{eq-cwt: variation of lambda_1}
	\end{align}
	If (\ref{eq-cwt: variation of lambda_1}) is zero for all $h$, then $\Delta(\varphi^2)=0\Rightarrow \varphi$ is constant, which implies that $\Sigma$ is a flat torus. Hence (\ref{eq-cwt: variation of lambda_1}) is nonzero for some $h$. We may change the sign of $h$ to make $d\lambda_1/dt$ positive. The perturbed metric $g(\delta)$ and eigenfunction $\varphi_\delta$ satisfies the theorem statement for small $\delta$.
\end{proof}

\subsection{Bi-Lipschitz equivalence with the round sphere}\label{subsec:bi-lipschitz}

\vspace{9pt}

Since $(\Sigma,\tilde g)$ is positively curved, its geometry is controlled by the homogeneous diameter-area ratio. In this subsection we bound the diameter of $\tilde g$ (recall that the area is normalized to be 1), and show that $\tilde g$ is bi-Lipschitz equivalent to the round sphere, then derive Theorem \ref{thm-intro:uniform-bi-Holder} as a corollary. The generic constants in this subsection depend only on $\beta$ and $|\Sigma|_g/\diam(\Sigma,g)^2$, and may vary from line to line.

\begin{theorem}\label{thm-cwt:diam_bound}
	$\diam(\Sigma,\tilde g)\leq C$, where $C$ depends on $\beta$ and the lower bound of $\frac{|\Sigma|}{\diam(\Sigma)^2}$.
\end{theorem}
\begin{proof}
	By (\ref{eq-cwt:K>e^-2u+du^2}) we have $\int|\D u|^2\,dA\leq 4\pi\beta^{-1}$. Denote $\bar u=\frac1{|\Sigma|}\int_\Sigma u\,dA$, by Theorem \ref{thm-intro:isop_ineq} and \ref{thm-aux:M-T-closed} we have
	\[\int_\Sigma e^{pu}\,dA\leq C|\Sigma|e^{p\bar u},\quad \int_\Sigma e^{-2u}\,dA\leq C|\Sigma|e^{-2\bar u}.\]
	In addition with Cauchy-Schwarz inequality $\int e^{2u}\cdot\int e^{-2u}\geq|\Sigma|^2$ and the normalization $\int e^{2u}\,dA=1$, we obtain $\int e^{pu}\,dA\leq C(\beta,p)|\Sigma|^{1-p/2}$ for all $p>1$.
	
	Denote $D=\diam(\Sigma,g)$, $A=|\Sigma|_g$. Let $x,y\in\Sigma$ such that $d_g(x,y)=D$. By coarea formula, there exists $s\in[\frac14D,\frac34D]$ such that $\gamma=\p B(x,s)$ is piecewise smooth and satisfies $|\gamma|\leq 4A/D$ and $\int_\gamma e^{2u}\leq 4/D$. Then Cauchy-Schwarz inequality gives $|\gamma|_{\tilde g}=\int_\gamma e^u\,dl\leq 4A^{1/2}/D$. Let $\gamma_1$ be the unique connected component of $\gamma$ that separates $x$ and $y$ (whose existence follows from $\pi_1(\Sigma)=0$), and denote $\Sigma\setminus\gamma_1=U\cup U'$. Assume $x\in U$, $y\in U'$. Theorem \ref{thm-intro:volume_comparison} gives $|U|\geq|B(x,s)|\geq CD^2$ and $|U'|\geq|B(y,\frac14D)|\geq CD^2$. Then any domain $\Omega\subset\subset U$ satisfies $|\Sigma\setminus\Omega|\geq|U'|\geq CD^2\geq C\frac{D^2}{A}\,|\Omega|$, therefore $\ID(U)\geq C\IN(\Sigma)$. Analogously, $\ID(U')\geq C\IN(\Sigma)$.
	
	The $\tilde g$-distance function $\tilde d(y)=d_{(U,\tilde g)}(y,\gamma_1)$ is Lipschitz, vanishing on $\p U=\gamma_1$, and satisfies $|\D\tilde d|_g=e^u$ almost everywhere. The $L^{n/(n-1)}$ Sobolev inequality gives
	\[||\tilde d||_{L^2(U,g)}\leq \ID(U)^{-1/2}||e^u||_{L^1(U,g)}.\]
	Following the iteration argument in \cite[Theorem 7.10]{Gilbarg-Trudinger_2001}, a uniform $L^{n/(n-1)}$ Sobolev inequality gives rise to a uniform $L^\infty$ Sobolev inequality
	\[||\tilde d||_{L^\infty}
		\leq C(p)\ID(U)^{-1/2}\,|\Sigma|^{1/2-1/p}\,||e^u||_{L^p(U,g)}
		\leq C,\qquad\forall p>2.\]
	The same conclusion holds for $\tilde d'=d_{(U',\tilde g)}(-,\gamma_1)$. Finally, the theorem follows from $\diam(\Sigma,\tilde g)\leq\max(\tilde d)+\max(\tilde d')+|\gamma_1|_{\tilde g}$.
\end{proof}

Next, we construct bi-Lipschitz maps from $(\Sigma,\tilde g)$ to the round sphere. By Weyl's embedding theorem and Lemma \ref{lemma-cwt: perturbation}, we assume that $(\Sigma,\tilde g)$ is the boundary of a convex body in $\RR^3$. We make the metric $\tilde g$ implicit when there is no ambiguity. The following lemma is elementary.

\begin{lemma}\label{lemma-cwt:biLip_from_inradius}
	Let $\Sigma\subset\RR^3$ be a convex surface, satisfying $d_1\leq|x|\leq d_2$ for any $x\in\Sigma$. Then $\Phi:x\mapsto\frac{x}{|x|}$ is a bi-Lipschitz map to the round sphere, with Lipschitz norms $||\Phi||_{Lip}\leq 1/d_1$, $||\Phi^{-1}||_{Lip}\leq d_2^2/d_1$.
\end{lemma}
\begin{proof}
	We compute
	\[d\Phi_x(v)=\frac{|x|^2v-(x\cdot v)x}{|x|^3},\quad v\in T_x\Sigma,\]
	and
	\[|d\Phi_x(v)|^2=\frac{|x|^2-|x\cdot v|^2}{|x|^4},\quad|v|=1.\]
	Thus $|d\Phi|\leq\frac{1}{|x|}\leq\frac1{d_1}$. To control $\Phi^{-1}$, we note that $x\cdot N\geq d_1$ by convexity. Hence for any $|v|=1$ we have
	\[|d\Phi_x(v)|^2=\frac{(x\cdot N)^2}{|x|^4}\geq\frac{d_1^2}{d_2^4}.\]
\end{proof}

\begin{prop}\label{prop-cwt:biLip_from_diameter}
	Let $\Sigma\subset\RR^3$ be a convex surface with unit area and (intrinsic) diameter $\leq D$. Then there exists a smooth bi-Lipschitz map to the round sphere $\Phi:\Sigma\to(S^2,g_0)$, with $||\Phi||_{Lip},||\Phi^{-1}||_{Lip}$ bounded in terms of $D$.
\end{prop}
\begin{proof}
	Let $W$ be the extrinsic width of $\Sigma$, defined as
	\begin{equation}
		W:=\min\big\{d: \text{we can embed }\Sigma\text{ into the region }\{-d/2\leq z\leq d/2\}\subset\RR^3\big\},
	\end{equation}
	where the embedding is isometric. First we consider the case that $W\geq\frac1{4\pi D}$. Let $S$ be a largest inscribed sphere of $\Sigma$, and $r$ be the radius of $S$. The width-inradius inequality gives
	\[\frac D2\geq r\geq \frac W4.\]
	A proof of the two-dimensional analogue can be found in \cite[p.215]{Jaglom-Boltjanskii}, while the proof there can be directly generalized to three dimensions. Now Lemma \ref{lemma-cwt:biLip_from_inradius} gives the desired bi-Lipschitz map.
	
	Next we assume $W\leq\frac1{4\pi D}$. We will replace the nearly flat top and bottom faces of $\Sigma$ by two cones, thus increasing the inradius. We may assume that $\Sigma$ is already embedded in the region $\{-W/2\leq z\leq W/2\}\subset\RR^3$. Let $Q$ be the projection image of $\Sigma$ onto the $xy$-plane, which is a convex planar domain. Thus $\Sigma\subset Q\times[-W/2,W/2]$. We have
	\[1=|\Sigma|\leq\big|\p(Q\times[-W/2,W/2])\big|\leq 2|Q|+2W|\p Q|\leq 2|Q|+2\pi DW\]
	Thus
	\[|Q|\geq\frac14.\]
	Let $r_0$ be the inradius of $Q$ (i.e. the radius of the largest disk contained in $Q$). From \cite{Scott-Awyong_2000} we know $r_0\geq \frac{|Q|}{2\diam(Q)}\geq \frac1{8D}$. After a translation, let
	\[B:=B(0,r_0)\subset Q\subset\RR^2\]
	be an inscribed disk in $Q$. Let $p^\pm=(0,0,\pm 1)$. Let $\Lambda^+$ be the convex cone obtained by taking the union of all tangent line segments of $\Sigma$ emanating from $p^+$. We claim that the normal unit vector $N_\Lambda$ of $\Lambda^+$ satisfies $|N_\Lambda\cdot\p_z|\geq\frac{r_0}{\sqrt{r_0^2+4}}$. To show this, let $P$ be any plane tangent to $\Lambda^+$. The intersection line $l=P\cap\{z=-1\}$ satisfies $d(l,p^-)\geq r_0$ by convexity. Our claim then follows from a simple computation. We can mollify the vertex of $\Lambda^+$ so that it is smooth and satisfies $|N_\Lambda\cdot\p_z|\geq\frac12\frac{r_0}{\sqrt{r_0^2+4}}$ everywhere.
	
	Let $\gamma^+$ be the set at which $\Lambda^+$ is tangent to $\Sigma$. Now $\gamma^+$ encloses a region $\Omega^+$ in $\Sigma$. Let $Q^+$ be the projection image of $\Omega^+$ onto the $xy$-plane. Consider the orthogonal projection map $F_1^+:\Omega^+\to Q^+$, $F_2^+:\Lambda^+\to Q^+$ (which are both bijective). It is not hard to prove the following general fact: if $G$ is the projection map from a surface $S$ to the $xy$-plane, then $||G^{-1}||_{Lip}\leq\sup_S|N_S\cdot\p_z|^{-1}$, where $N_S$ is the unit normal vector of $S$. By our claim above, the map $F^+=(F_2^+)^{-1}\circ F_1^+: \Omega^+\to\Lambda^+$ has controlled bi-Lipschitz norms. Analogously, we can define a bi-Lipschitz map $F^-:\Omega^-\to\Lambda^-$
	
	Let $\Omega^m=\Sigma\setminus(\Omega^+\cup\Omega^-)$.
	Combining $F^+$ with $F^-$ we obtain a $C^{1,1}$ map $F$ from $\Sigma$ to $\Sigma'=\Omega^m\cup\Lambda^+\cup\Lambda^-$, with $||F||_{Lip}\leq C$, $||F^{-1}||_{Lip}\leq C$. It is not hard to see $\diam(\Sigma')<\diam(\Sigma)+4$ and $|x|>\min(1/4,r_0/4)$ on $\Sigma'$. From Lemma \ref{lemma-cwt:biLip_from_inradius} we obtain a $C^{1,1}$ map $\Phi_0:\Sigma\to S^2$ with controlled bi-Lipschitz norms. Finally, we mollify $\Phi_0$ to obtain a smooth map $\Phi$ that $C^1$-approximates $\Phi_0$ and is therefore bi-Lipschitz.
\end{proof}

\begin{proof}[Proof of Theorem \ref{thm-intro:uniform-bi-Holder}] {\ }
	
	We have obtained three metrics on $\Sigma$: the original metric $g$, the conformal metric $\tilde g=e^{2u}g$, and the pulled-back round metric $\Phi^*g_0$ which satisfies $C_1^{-1}\tilde g\leq \Phi^*g_0\leq C_1\tilde g$ by Proposition \ref{prop-cwt:biLip_from_diameter}. The constant $C_1$ depends on $\diam(\tilde g)$, hence depends on the lower bound $A_0$ of $|\Sigma|_g$, under the normalizing assumptions of the theorem.
	
	Consider the metric $g_1=e^{-2u_1}g_0$ where $u_1=u\circ\Phi^{-1}$. We have $C_1^{-1}(\Phi^{-1})^*g\leq g_1\leq C_1(\Phi^{-1})^*g$ hence $g_1$ is bi-Lipschitz equivalent to $g$. It remains to prove that $g_1$ and $g_0$ are uniformly bi-H\"older equivalent. This directly follows from \cite[Proposition 6.3]{Aldana-Carron-Tapie_2021} (whose notations can be found at the beginning of Section 4A). By (\ref{eq-cwt:K>e^-2u+du^2}) we have $\int|\tilde\D u|^2\,d\tilde A\leq 4\pi\beta^{-1}$, hence $\int|\D_0 u_1|^2\,dA_0\leq 4\pi\beta^{-1}C_1^2$ is uniformly bounded. This verifies the assumption in the cited theorem. Moreover,
	\[C_1^{-1}A_0\leq C_1^{-1}|\Sigma|_g\leq|\Sigma|_{g_1}\leq C_1|\Sigma|_g\leq C\]
	by our assumptions and volume comparison, hence \cite[Proposition 6.3]{Aldana-Carron-Tapie_2021} gives uniform bi-H\"older equivalence between $g_0$ and $g_1$.
\end{proof}

\begin{remark}\label{rmk-cwt:discussion}
	Underlying the bi-H\"older equivalence result is the deep theory of strong $A_\infty$ weights, initially introduced by David and Semmes \cite{David-Semmes_1990} and developed recently by Aldana, Carron and Tapie \cite{Aldana-Carron-Tapie_2021} for closed manifolds. A density function $e^{nf}$ on an $n$-manifold $M$ is called a strong $A_\infty$ weight if the metric $g'=e^{2f}g$ satisfies
	
	(1) volume-doubling condition: $|B_g(x,2r)|_{g'}\leq C|B_g(x,r)|_{g'}$,
	
	(2) volume-distance compatibility: $C^{-1}d_{g'}(x,y)^n\leq\big|B_g(x,d_g(x,y))\big|_{g'}\leq Cd_{g'}(x,y)^n$, \\
	as defined in \cite[Definition 4.1]{Aldana-Carron-Tapie_2021}. The second condition seems somewhat unnatural in the sense that $g$-metric balls appear in the central term of the inequality. However, it is true that $g'$ satisfies Euclidean volume growth, i.e.
	\[(C')^{-1}d_{g'}(x,y)^n\leq\big|B_{g'}(x,d_{g'}(x,y))\big|_{g'}\leq (C')d_{g'}(x,y)^n.\]
	Moreover, it is shown by David and Semmes \cite{David-Semmes_1990} that $g'$ satisfies a uniform isoperimetric inequality and Sobolev inequality.
	
	It is an important fact that functions in the critical Sobolev space $W^{k,n/k}$ are strong $A_\infty$ weights, see \cite{Bonk-Heinonen-Saksman_2004} for the case of $\RR^n$ and \cite{Aldana-Carron-Tapie_2021} for the case of closed manifolds. In our situation $u_1\in W^{1,2}$, therefore $e^{-2u_1}$ is a strong $A_\infty$ weight for $g_0$. David and Semmes' isoperimetric inequality thus exactly corresponds to Theorem \ref{thm-intro:isop_ineq}. However, this does not constitute an independent proof of theorem \ref{thm-intro:isop_ineq} since the uniform control on $\tilde g$ is based on them.
\end{remark}

\subsection{The case of almost rigidity}\label{subsec:rig}

\vspace{9pt}

In this subsection we prove Theorem \ref{thm-intro:almost_rigidity}. The H\"older constants originating from strong $A_\infty$ weights are not sharp, so the proof requires more efforts. By assumption, the first eigenfunction $\varphi>0$ satisfies $\Delta \varphi\leq(\beta K-\beta)\varphi$. From (\ref{eq-intro:eigenvalue_condition_main}) one obtains $|\Sigma|_g\leq 4\pi$. The same conformal transformation in subsection \ref{subsec:preparation} now gives
\begin{equation}\label{eq-rig:K>e^-2u+du^2}
	\tilde K\geq e^{-2u}+\beta|\tilde\D u|^2.
\end{equation}
This time we adopt the normalization $|\Sigma|_{\tilde g}=\int_\Sigma e^{2u}\,dA=4\pi$ for the purpose of matching the areas of $g$ and $\tilde g$. Integrating (\ref{eq-rig:K>e^-2u+du^2}) against the volume form of $\tilde g$ and using the condition $|\Sigma|_g\geq 4\pi-\delta$, we obtain $\int_\Sigma|\tilde\D u|^2\,d\tilde A\leq\beta^{-1}\delta$. \\

\textbf{Notation.} In this subsection, we use $o_\delta(1)$ to denote universal constants that depends only on $\beta,\delta$ and converges to zero when $\delta\to0$. For example, $4\pi-o_\delta(1)\leq\int e^{-2u}\,d\tilde A\leq4\pi$. The generic constants $C$ depend only on $\beta$ unless explicitly indicated.

\begin{lemma}\label{lemma-rig:bar_u}
	Denote $\bar u=\frac1{4\pi}\int_\Sigma u\,d\tilde A$. Then $|\bar u|=o_\delta(1)$. Moreover, for any fixed $\eta>0$ we have $|\{u\geq\eta\}|=o_\delta(1)$.
\end{lemma}
\begin{proof}
	A simple estimate gives
	\[4\pi\geq\int_\Sigma e^{-2u}\,d\tilde A\geq\int_\Sigma(1-2u)\,d\tilde A\ \Rightarrow\ \bar u\geq0.\]
	The other direction requires finite energy of $u$. By Theorem \ref{thm-intro:bonnet_myers}, \ref{thm-cwt:diam_bound} and Burago-Zalgaller's inequality, $\IN(\Sigma,\tilde g)$ has a uniform lower bound. Theorem \ref{thm-aux:M-T-closed} gives
	\[\int_\Sigma\exp\Big[\frac{(u-\bar u)^2}{C\delta}\Big]\,d\tilde A\leq C.\]
	Therefore,
	\[\int_\Sigma\exp\big(\delta^{-1/2}|u-\bar u|\big)\,d\tilde A
	\leq \int_\Sigma\exp\Big[\frac{(u-\bar u)^2}{C\delta}+\frac14C\Big]\,d\tilde A\leq C.\]
	Set $S_1=\{|u-\bar u|\leq\delta^{1/4}\}$ and $S_2=\Sigma\setminus S_1$. We have
	\[\begin{aligned}
		4\pi-\delta &\leq \int_\Sigma e^{-2u}\,d\tilde A
		= \int_{S_1}e^{-2u}\,d\tilde A+\int_{S^2}e^{-2u}\,d\tilde A \\
		&\leq 4\pi e^{-2\bar u+2\delta^{1/4}}
			+ e^{-2\bar u}\int_{S_2}e^{|u-\bar u|/\sqrt\delta}\,d\tilde A
			\cdot \exp\big[(2-\delta^{-1/2})\delta^{1/4}\big] \\
		&= e^{-2\bar u}(4\pi+o_\delta(1))
	\end{aligned}\]
	Hence $\bar u\geq -o_\delta(1)$, and this shows $|\bar u|=o_\delta(1)$. Finally, for any fixed $\eta>0$ we have by Chebyshev's inequality
	\begin{equation}\label{eq-rig:small_measure}
		\Big|\{|u|\geq\eta\}\Big|\cdot e^{\eta/\sqrt\delta}\leq e^{|\bar u|/\sqrt\delta}\int_\Sigma\exp\big(\delta^{-1/2}|u-\bar u|\big)\,d\tilde A.
	\end{equation}
	This shows $|\{|u|\geq\eta\}|\to0$ when $\delta\to0$.
\end{proof}

\begin{prop}\label{prop-rig:bi-Lip-almost_rig}
	Embed $(\Sigma,\tilde g)$ as the boundary of a convex body in $\RR^3$. After an appropriate translation, the map $\Phi:x\mapsto\frac{x}{|x|}$ satisfies $||\Phi||_{Lip}\leq 1+o_\delta(1)$ and $||\Phi^{-1}||_{Lip}\leq 1+o_\delta(1)$.
\end{prop}
\begin{proof}
	All the integrals in this proof are with respect to $\tilde g$, so we omit the area form for brevity. Fix $\eta\ll1$. By Lemma \ref{lemma-rig:bar_u}, for $\delta$ sufficiently small we have
	\[|\{K\leq 1-\eta\}|\leq\eta.\]
	Denote $S_1=\{K\leq1-\eta\}$, $S_2=\{K\geq1+\sqrt\eta\}$, $S_3=\Sigma\setminus(S_1\cup S_2)$. By Gauss-Bonnet formula we have
	\[\begin{aligned}
		4\pi &= \int_{S_1}K+\int_{S_2}K+\int_{S_3}K
		\geq \int_{S_2}K+(4\pi-\eta-|S_2|)(1-\eta) \\
		&\geq \int_{S_2}K+(4\pi-\eta-\frac{\int_{S_2}K}{1+\sqrt\eta})(1-\eta).
	\end{aligned}\]
	Hence $\int_{S_2}K\leq C\sqrt\eta$.
	
	We first translate $\Sigma$ so that the origin is contained in the interior. Let $N$ be the outer unit normal of $\Sigma$, and $H$ be the mean curvature. Let $V$ be the volume enclosed by $\Sigma$. Minkowski' inequality $\int_\Sigma H\geq4\sqrt{\pi|\Sigma|_{\tilde g}}$ and integration formula $\int_\Sigma K(x\cdot N)=\frac12\int_\Sigma H$ (see \cite{Hsiung_1954}; $x$ is the position vector) yields
	\[\begin{aligned}
		4\pi &\leq \frac12\int_\Sigma H
		= \int_\Sigma K(x\cdot N)
		= \int_{S_2}K(x\cdot N)+\int_{\Sigma\setminus S_2}K(x\cdot N) \\
		&\leq C\sqrt\eta\cdot\diam(\Sigma,\tilde g)
			+ (1+\sqrt\eta)\int_\Sigma(x\cdot N) \\
		&= C\sqrt\eta+3(1+\sqrt\eta)V,
	\end{aligned}\]
	where the last line follows from divergence formula. Hence
	\begin{equation}\label{eq-s2s1:upper_bound_V}
		V\geq \frac43\pi-C\sqrt\eta.
	\end{equation}
	Let $\rho$ be the inradius of $\Sigma$ (i.e. the radius of the largest ball enclosed by $\Sigma$). By a Bonnesen-type isoperimetric inequality, see formula (114) in \cite{Osserman_1979}, we have
	\[\Big(\frac{|\Sigma|}{4\pi\rho^2}\Big)^{3/2}-\Big(\frac{V}{\frac43\pi\rho^3}\Big)
	\geq \Big[\big(\frac{|\Sigma|}{4\pi\rho^2}\big)^{1/2}-1\Big]^3
	\ \Rightarrow\ (1-\rho)^3\leq 1-\frac{V}{\frac43\pi}.\]
	Combined with (\ref{eq-s2s1:upper_bound_V}) this gives
	\[\rho\geq 1-C\eta^{1/6}.\]
	Translate $\Sigma$ so that $\min_{x\in\Sigma}|x|\geq 1-C\eta^{1/6}$. It remains to show $\max_{x\in\Sigma}|x|\leq 1+o_\eta(1)$. Suppose $x\in\Sigma$ with $|x|\geq R$. Let $\Omega_1$ be the solid ball with radius $\rho$, and $\Omega_2$ be the cone with vertex $x$ and is tangent to $\Omega_1$. Their union $\Omega=\Omega_1\cup\Omega_2$ is a convex body. The orthogonal projection $\Sigma\to\p\Omega$ is 1-Lipschitz and does not increase area. Hence
	\[\begin{aligned}
		4\pi &= |\Sigma|\geq|\p\Omega|
		= \rho^2 \Big[2\pi(1+\cos\th)+\pi\sin\th\tan\th\Big] \\
		&\geq (1-C\eta^{1/6})^2 \,4\pi\, (1+\frac{\sin^4(\th/2)}{\cos\th}),
	\end{aligned}\]
	where $\th=\arccos(\rho/R)$. Therefore $R\leq1+o_\eta(1)$ when $\eta\to0$. Now the proposition follows from Lemma \ref{lemma-cwt:biLip_from_inradius}.
\end{proof}

\begin{proof}[Proof of Theorem \ref{thm-intro:almost_rigidity}] {\ }
	
	Similar to the proof of Theorem \ref{thm-intro:uniform-bi-Holder}, we let $u_1=u\circ\Phi^{-1}$ and $g_1=e^{-2u_1}g_0$. Proposition \ref{prop-rig:bi-Lip-almost_rig} implies that $g_1$ is $\big(1+o_\delta(1)\big)$-bi-Lipschitz equivalent to $g$. The remaining work is showing $|d_{g_1}(x,y)-d_{g_0}(x,y)|=o_\delta(1)$ uniformly for all $x\ne y$.
	
	Note that $\int_\Sigma|\D_0 u_1|^2\,dA_0\leq C\delta$ and $|\int e^{-2u_1}\,dA_0-4\pi|=o_\delta(1)$ by Lipschitz equivalence between $\tilde g$ and $g_0$. Arguing as in Theorem \ref{thm-cwt:diam_bound}, we have $\int_\Sigma e^{-pu_1}\,dA_0\leq C(\beta,p)$. The function $d(y)=d_{g_1}(x,y)$ is Lipschitz and satisfies $|\D_0 d|_{g_0}=e^{-u_1}$. By Morrey's inequality we have
	\begin{equation}\label{eq-rig:morrey}
		d_{g_1}(x,y)\leq Cd_{g_0}(x,y)^\alpha,\quad \forall\,\text{fixed }0<\alpha<1.
	\end{equation}
	
	
	We first show the more complicated direction $d_{g_1}(x,y)\geq d_{g_0}(x,y)-o_\delta(1)$. Suppose that $d_{g_1}(x,y)\leq d_{g_0}(x,y)-s$ for some $s>0$ independent of $\delta$ and for some $x\ne y$. Let $\gamma$ be a shortest $g_1$-geodesic from $x$ to $y$. Denote $d_0=d_{g_0}(x,y)$, $a=\frac{d_0-s}{d_0-s/2}$. We claim that we can assume $d_0\leq\frac23\pi$. Suppose $d_0\geq\frac23\pi$, we embed $\Sigma$ as the unit sphere so that $x,y$ has the same longitude and opposite latitude ($=\pm d_0/2$). Choose any $z\in\gamma$ that lies on the equator, thus $d_{g_0}(x,z)=d_{g_0}(y,z)\leq\pi-d_0/2\leq \frac23\pi$. Moreover,
	\[d_{g_1}(x,z)+d_{g_1}(y,z)=d_{g_1}(x,y)\leq d_{g_0}(x,z)+d_{g_0}(y,z)-s.\]
	Without loss of generality we assume $d_{g_1}(x,z)\leq d_{g_0}(x,z)-s/2$. Then we can proceed with the pair of points $(x,z)$ and new constant $s/2$.
	
	Let $U=\{z\in\gamma:e^{-u_1}(z)>a\}$, $V=\gamma\setminus U$, hence $|U|_{g_0}\leq d_0-s/2$. We would like to apply the Moser-Trudinger inequality to $e^{-2u_1}$ with Dirichlet condition on $V$. Therefore, we need a corresponding isoperimetric inequality:
	
	
	\begin{lemma}\label{lemma-rig:IN_bound}
		Any domain $\Omega\subset\Sigma$ with $\bar\Omega\cap V=\emptyset$ satisfies $|\p\Omega|_{g_0}^2\geq \xi|\Omega|_{g_0}$ for some constant $\xi$ depending only on $\beta$ and the lower bound of $s$.
	\end{lemma}
	We postpone its proof to the end of this section. Let $b=\min_V(u_1)$, so $e^{-2b}\leq a^2\leq\big(\frac{\pi-s}{\pi-s/2}\big)^2$. Apply Theorem \ref{thm-aux:M-T-Dirichlet} to $v=\max\big\{\frac{b-u_1}{\sqrt\delta},0\big\}$, which satisfies $\int_\Sigma|\D_0 v|^2\,dA_0\leq C$:
	\[\int_\Sigma\exp\big(\frac{b-u_1}{\sqrt\delta}\big)\,dA_0
		\leq \int_\Sigma e^v\,dA_0
		\leq \int_\Sigma \exp\Big(\frac{\xi v^2}{\int|\D_0 v|^2}+\frac1{4\xi}\int_\Sigma|\D_0 v|^2\,dA_0\Big)
		\leq C,\]
	where $\xi$ is the uniform lower bound of $\ID(\Sigma\setminus V,g_0)$. Let $S_1=\{u_1\geq b-\delta^{1/4}\}$, $S_2=\Sigma\setminus S_1$. We have
	\[\begin{aligned}
		\int_\Sigma e^{-2u_1}\,dA_0 &\leq e^{-2b}
			\Big( |S_1|e^{2\delta^{1/4}} + \int_{S_2}e^{2(b-u_1)}\,dA_0 \Big) \\
		&\leq \big(\frac{\pi-s}{\pi-s/2}\big)^2
			\Big(4\pi e^{2\delta^{1/4}}
			+\exp\big[(2-\delta^{-1/2})\delta^{1/4}\big]\int_{S_2}e^{(b-u_1)/\sqrt\delta}\,dA_0\Big) \\
		&\leq \big(\frac{\pi-s}{\pi-s/2}\big)^24\pi(1+o_\delta(1)).
	\end{aligned}\]
	This contradicts $|\Sigma|_{g_1}\geq 4\pi-o_\delta(1)$ when $\delta$ is sufficiently small. Therefore $d_{g_1}(x,y)\geq d_{g_0}(x,y)-o_\delta(1)$.
	
	Next we prove that $d_{g_1}(x,y)\leq d_{g_0}(x,y)+o_\delta(1)$. Assume otherwise that $d_{g_1}(x,y)\geq d_{g_0}(x,y)+s$ for some $s>0$ independent of $\delta$. By (\ref{eq-rig:morrey}) this implies $d_{g_0}(x,y)\geq C(\beta,s)$. Note that $\int_\Sigma|\D_0 u_1|\,dA_0\leq C\sqrt\delta$. By coarea formula, there exists a regular value $s\in[\bar u-2\delta^{1/4},\bar u-\delta^{1/4}]$ such that $\big|\{u=s\}\big|\leq C\delta^{1/4}$. The argument in (\ref{eq-rig:small_measure}) implies $\big|\{u\leq s\}\big|=o_\delta(1)$. Combining these two conditions, any geodesic segment contained in $\{u\leq s\}$ must have length $\leq o_\delta(1)$. Let $\gamma$ be the shortest $g_0$-geodesic joining $x$ and $y$. Consider the parallel curves $\gamma_t=\exp_\gamma(t\nu)$, where $\nu$ is the nuit normal vector field of $\gamma$ with any orientation. By coarea formula, there exists a value $t_0\in[-C\delta^{1/4},C\delta^{1/4}]$ such that $\gamma_{t_0}\cap\{u=s\}=\emptyset$. Hence $\gamma_{t_0}\subset\{u\geq s\}$ when $\delta$ is sufficiently small. Let $x',y'$ be the endpoints of $\gamma_{t_0}$. We now have
	\[\begin{aligned}
		d_{g_1}(x,y) &\leq d_{g_1}(x,x')+d_{g_1}(y,y')+|\gamma_{t_0}|_{g_1} \\
		&\leq 2C\delta^{\alpha/4}+e^{-2s}|\gamma_{t_0}|_{g_0} \\
		&\leq d_{g_0}(x,y)+o_\delta(1),
	\end{aligned}\]
	which is a contradiction. This completes the proof.
\end{proof}

\begin{proof}[Proof of Lemma \ref{lemma-rig:IN_bound}] {\ }
We may assume that $\Omega$ is connected, as seen in the proof of Theorem \ref{thm-intro:isop_ineq}. Let $\{\phi,\th\}$ be the spherical coordinate system. After a rotation we may assume $x,y$ have coordinates $(-d_0/2,0), (d_0/2,0)$. Let $P_\phi:\Sigma\to[-\pi/2,\pi/2]$ map a point to its latitude, and define $I=P_\phi(V)\cap[-d_0/2,d_0/2]$. Since $P_z$ is a surjective from $\gamma$ to $[-d_0/2,d_0/2]$, we have $|I|\geq d_0-|P_\phi(U)|\geq d_0-|U|_{g_0}\geq s$. Let $J=P_\phi(\p\Omega)$. If $J\supset I$, then
\[|\p\Omega|_{g_0}\geq|J|\geq|I|\geq s\geq\frac{s}{\sqrt{4\pi}}|\Omega|_{g_0}^{1/2}.\]
Now suppose $J\nsupseteq I$, then there exists $\phi_1\in I$ with $\p\Omega\cap\{\phi=\phi_1\}=\emptyset$. Therefore, we have either $\{\phi=\phi_1\}\subset\Omega$ or $\{\phi=\phi_1\}\cap\Omega=\emptyset$. The first case cannot happen since we have assumed $\bar\Omega\cap V=\emptyset$, hence the second case must hold. Now, $\Omega$ is contained in either the spherical cap $\{\phi>\phi_1\}$ or $\{\phi<\phi_1\}$. Since $\phi_1\subset[-\pi/3,\pi/3]$, the conclusion follows from the usual spherical isoperimetric inequality.
\end{proof}

\appendix

\section{Two Proofs of weak Bonnet-Myers' Theorem}\label{sec:bonnet_myers}

In this section we present two proofs of Theorem \ref{thm-intro:bonnet_myers}, using the methods of weighted geodesics and weighted $\mu$-bubbles. \\

Let $u>0$ satisfies
\begin{equation}
	\Delta u^\beta\leq(\beta K-\lambda)u^\beta,
\end{equation}
equivalently,
\begin{equation}\label{eq-appA:eq_of_u}
	\Delta u\leq(K-\lambda\beta^{-1})u+(1-\beta)u^{-1}|\D u|^2.
\end{equation}
Note that the function $u$ defined here is different than in Section \ref{sec:cwt}.

\subsection{Proof using weighted geodesics}

\vspace{9pt}

The argument here is similar to the ones in Schoen-Yau \cite{Schoen-Yau_1983} and Shen-Ye \cite{Shen-Ye_1996}.

Let $p,q\in\Sigma$ be two points with the largest distance. If $\Sigma$ is non-compact, then choose $p,q$ with large enough distance to obtain a contradiction below. Let $\gamma:[0,L]\to\Sigma$, $\gamma(0)=p$, $\gamma(L)=q$, be a minimizer of the weighted length functional $\int_\gamma u\,dl$. Parametrize $\gamma$ with unit speed. The second variational formula gives the following inequality:
\begin{equation}\label{eq-appA:length_2nd_var}
	0\leq\int_\gamma\Big[u(\varphi')^2+(\Delta u-u'')\varphi^2-Ku\varphi^2-u^{-1}(u_N)^2\varphi^2\Big]\,dl
\end{equation}
for any function $\varphi:\varphi(0)=\varphi(L)=0$, where we denote $f':=\p f/\p\gamma'$ for a function $f$, and denote $u_N:=\p u/\p N$. Substituting $\varphi=u^{-1/2}\psi$ into (\ref{eq-appA:length_2nd_var}) and using equation (\ref{eq-appA:eq_of_u}), we obtain
\[\begin{aligned}
	0 &\leq \int_\gamma\Big[u\big(-\frac12u^{-3/2}u'\psi+u^{-1/2}\psi'\big)^2+u'(-u^{-2}u'\psi^2+2u^{-1}\psi\psi') \\
	&\qquad\qquad -\lambda\beta^{-1}\psi^2+(1-\beta)u^{-2}(u')^2\Big]\,dl \\
	&\leq \int_\gamma\Big[(\frac14-\beta)u^{-2}(u')^2\psi^2+(\psi')^2+u^{-1}u'\psi\psi'-\lambda\beta^{-1}\psi^2\Big]\,dl \\
	&\leq \int_\gamma\Big[\big(1+\frac14(\beta-\frac14)^{-1}\big)(\psi')^2-\lambda\beta^{-1}\psi^2\Big]\,dl.
\end{aligned}\]
Theorem \ref{thm-intro:bonnet_myers} follows by letting $\psi(t)=\sin(\pi t/L)$.

\subsection{Proof using weighted \texorpdfstring{$\mu$}{μ}-bubbles}

\vspace{9pt}

A $\mu$-bubble is a hypersurface with prescribed mean curvature $H=h$, where $h$ is a given function on the ambient manifold. The variational characterization of $\mu$-bubbles is given by the energy functional (\ref{eq-appA:energy_for_mu_bubbles}). By choosing an appropriate function $h$, we can extract information about the ambient manifold from the stability inequality. The introduction here closely follows Chodosh-Li \cite{Chodosh-Li}.

Let $\Omega^+$, $\Omega^-$ be two disjoint domains in $\Sigma$. Let $h$ be a Lipschitz function (whose conditions will be determined later) on $\Sigma\setminus(\Omega^-\cup\Omega^+)$, such that $h|_{\p\Omega^\pm}=\pm\infty$. Let $\Omega^0$ be a domain containing $\Omega^+$ and disjoint from $\Omega^-$. (This set serves the role of renormalizing.) Consider the following functional acting on all open sets $\Omega$ with $\Omega\Delta\Omega^0\subset\subset\Sigma\setminus(\Omega^-\cup\Omega^+)$:
\begin{equation}\label{eq-appA:energy_for_mu_bubbles}
	E_u(\Omega):=\int_{\p\Omega}u\,dl-\int_M(\chi_\Omega-\chi_{\Omega^0})hu\,dA.
\end{equation}
A critical point of $E(\Omega)$ (or its boundary) is called a $\mu$-bubble. Since $h$ is infinite on $\p\Omega^\pm$, any $\mu$-bubble must lie between $\Omega^+$ and $\Omega^-$. Note that we have added a weight $u$ into the functional. The unweighted version is $E(\Omega)=|\p\Omega|-\int_M(\chi_\Omega-\chi_{\Omega^0})h\,dA$, whose critical point satisfies $H=h$. By geometric measure theory, a $\mu$-bubble is always a $C^{2,\alpha}$ hypersurface.

In \cite{Chodosh-Li} is was shown in detail that a global minimizer of (\ref{eq-appA:energy_for_mu_bubbles}) always exists. The first variation of (\ref{eq-appA:energy_for_mu_bubbles}) gives
\begin{equation}\label{eq-appA: solve1}
	\kappa=h-u^{-1}u_N
\end{equation}
for a $\mu$-bubble, where $\kappa$ is the geodesic curvature of $\p\Omega$. The second variation of (\ref{eq-appA:energy_for_mu_bubbles}) at a global minimizer gives the following stability inequality:
\begin{equation}\label{eq-appA:mu_bubble_2nd_var}
	\begin{aligned}
		0 &\leq \int_\gamma\Big[u(\varphi')^2-Ku\varphi^2-h^2u\varphi^2+hu_N\varphi^2-u^{-1}u_N^2\varphi^2 \\
		&\qquad\qquad +(\Delta u-u'')\varphi^2-h_Nu\varphi^2\Big]\,dl.
	\end{aligned}
\end{equation}
Testing (\ref{eq-appA:mu_bubble_2nd_var}) with $\varphi=u^{-1}$, we obtain
\[\begin{aligned}
	0 &\leq \int_\gamma\Big[-\beta u^{-3}(u')^2+\big(-h^2-\lambda\beta^{-1}+|\D h|\big)u^{-1}+hu^{-2}u_N-\beta u^{-3}u_N^2\Big]\,dl \\
	&\leq \int_\gamma\Big[(\frac1{4\beta}-1)h^2-\lambda\beta^{-1}+|\D h|\Big]u^{-1}\,dl.
\end{aligned}\]
For the last line we used
\begin{equation}\label{eq-appA:solve2}
	hu^{-2}u_N\leq\beta u^{-3}u_N^2+\frac1{4\beta}h^2u^{-1}.
\end{equation}
Therefore, we would reach a contradiction if $h$ satisfies
\begin{equation}\label{eq-appA:solve3}
	|\D h|<(1-\frac1{4\beta})h^2+\lambda\beta^{-1}.
\end{equation}
Suppose $\beta>\frac14$, the function
\[h(x)=\sqrt{\frac{\lambda}{\beta(1-\frac{1}{4\beta}-\epsilon)}}\cot\Big[\sqrt{\lambda\beta^{-1}(1-\frac1{4\beta}-\epsilon)}\,d(x,\Omega^+)\Big]=:C_1\cot\big[C_2\,d(x,\Omega^+)\big]\]
satisfies (\ref{eq-appA:solve3}). If $\diam(\Sigma)>\pi/C_2$, then we can let $p,q$ be points with maximal distance, and choose $\Omega^+=B_\epsilon(p)$, $\Omega^-=\{x\in\Sigma: d(x,\Omega^+)\geq\pi/C_2\}$. For these choices we obtain a contradiction with the stability inequality. This proves $\diam(M)\leq \pi/C_2$.

\section{Counterexamples for \texorpdfstring{$\frac14<\beta\leq\frac12$}{1/4<β≤1/2}}\label{sec:countereg}

\vspace{9pt}

In this section we discuss the sharpness of Theorem \ref{thm-intro:isop_ineq} and \ref{thm-intro:volume_comparison} for $\frac14<\beta\leq\frac12$ by finding concrete counterexamples. The main piece for the counterexample is
\[g=dr^2+r^{-p}d\th^2,\quad \varphi=r^q\qquad(p>0,r>0).\]
A short computation shows that $\Delta\varphi\leq\beta K\varphi$ is equivalent to
\begin{equation}\label{eq-intro: aux1}
	q(q-1)-pq+\beta p(p+1)\leq0.
\end{equation}
For $\beta>\frac14$, the largest achievable value of $p$ under (\ref{eq-intro: aux1}) is $p=\frac1{4\beta-1}$. Hence when $\beta\leq\frac12$ we can achieve $p\geq1$. Observe the following properties of $g$ when $p\geq1$: the end at $r\to0$ has infinite area, the perimeter of the end at $r\to\infty$ has the order $r^{-p}$, the end at $r\to\infty$ has finite area ($p>1$) or has area growth $\sim\log r$ ($p=1$). From $g$ we can build the following counterexamples:

\begin{example}\label{ex-countereg:volume}
	Truncate $g$ at a fixed large $r$ and glue a spherical cap with a properly chosen metric and $\varphi$, such that $\Delta\varphi\leq\beta K\varphi$ is satisfied after the gluing. In this way we obtain an non-compact surface that violates Theorem \ref{thm-intro:volume_comparison}.
\end{example}

\begin{example}\label{ex-countereg:collapsing}
	Truncate $g$ at a large radius $R$ and a fixed small radius $r_0$. Glue a spherical cap at $r_0$ and a half catenoid-shaped metric (see construction below) at $R$. Denote the resulting surface with boundary by $\Omega$. Taking the double of $\Omega$ we obtain a closed surface $\Sigma$. Then $\diam(\Sigma)=O(R)$ as $R\to\infty$, and the area of $\Omega$ has the order
	\[|\Omega|=\left\{\begin{aligned}
		& O(R^{1-p})\qquad(p<1), \\
		& O(\log R)\qquad(p=1), \\
		& O(1)\qquad(p>1).
	\end{aligned}\right.\]
	Hence,
	\[\Ch(\Sigma)\cdot\diam(\Sigma)\leq\frac{|\p\Omega|}{|\Omega|}\cdot\diam(\Sigma)=\left\{\begin{aligned}
		& O(R^{1-p})\qquad(p>1), \\
		& O(1/\log R)\qquad(p=1).
	\end{aligned}\right.\]
	Since $\Ch(\Sigma)\cdot\diam(\Sigma)$ is scale-invariant, by scaling this example we obtain closed surfaces with unit diameter and arbitratily small Cheeger's constant. This verifies the remark after Theorem \ref{thm-intro:isop_ineq}. When $\beta<\frac12$, for the largest possible $p$ we have $\IN(\Sigma)\leq\frac{|\p\Omega|^2}{|\Omega|}=O(R^{-2p})=O(R^{-\frac2{4\beta-1}})$ and $\frac{|\Sigma|}{\diam(\Sigma)^2}=O(R^{-2})$, hence the power in Theorem \ref{thm-intro:isop_ineq}(2) is almost sharp.
\end{example}

Below is the technical construction of Example \ref{ex-countereg:collapsing}. Example \ref{ex-countereg:volume} can be similarly constructed, and we omit the details. Assume $\beta<\frac12$ for simplicity. Define the follwing constants:
\[p=2-2\beta\ (>1),\quad
	c_1=\frac16,\quad c_2=\frac35,\quad c_3=40,\quad
	A=\big(\frac{c_3p}{r_0}\big)^2,\]
\[r_0=\big(c_2p\sqrt{1+c_3^2}\big)^{\frac{1}{p+1}},\quad
	r_1=r_0-\frac1{\sqrt A}\big(\pi-\tan^{-1}(c_3)-\cos^{-1}(c_2)\big).\]
It is mentioned above that $=\frac1{4\beta-1}$ is the largest achievable value of $p$. However, Example \ref{ex-countereg:collapsing} holds for all $p>1$, therefore we set $p=2-2\beta$ for simplicity of expressions. Consider the metric $g=dr^2+f^2(r)d\th^2$ and function $\varphi=\varphi(r)$ defined as follows:
\[f(r)=\left\{\begin{aligned}
	& \frac{p}{2c_1}R^{-p-2}(r-R-c_1R)^2+(1-\frac12c_1p)R^{-p}\qquad(R\leq r\leq R+c_1R) \\
	& r^{-p}\qquad(r_0\leq x\leq R) \\
	& r_0^{-p}\cos\big[\sqrt A(r-r_0)\big]-r_0^{-p-1}\frac{p}{\sqrt{A}}\sin\big[\sqrt A(r-r_0)\big]\qquad(r_1\leq r\leq r_0) \\
	& f(r_1^+)+r-r_1\qquad(r_1-f(r_1^+)\leq r\leq r_1)
\end{aligned}\right.\]
\[\varphi(r)=\left\{\begin{aligned}
	& -\frac1{2c_1R}(r-R-c_1R)^2+(1+\frac12c_1)R\qquad(R\leq r\leq R+c_1R) \\
	& r\qquad(r_0\leq r\leq R) \\
	& r_0\cosh\big[\frac{\sqrt{\beta A}}{2}(r-r_0)\big]+\frac2{\sqrt{\beta A}}\sinh\big[\frac{\sqrt{\beta A}}{2}(r-r_0)\big]\qquad(r_1\leq r\leq r_0) \\
	& \varphi(r_1^+)\qquad(r_1-f(r_1^+)\leq r\leq r_1)
\end{aligned}\right.\]

\begin{lemma}
	$f(r)$ is a $C^{1,1}$ function on $[r_1-f(r_1^+),R+c_1R]$. $\varphi(r)$ is positive and $C^{1,1}$ except at $r=r_1$, where one has $\varphi'(r_1^-)=0>\varphi'(r_1^+)$. For $R$ sufficiently large, in each interval $(r_1-f(r_1^+),r_1)$, $(r_1,r_0)$, $(r_0,R)$, $(R,R+c_1R)$ we have
	\begin{equation}\label{eq-countereg:supersol}
		\varphi''+\frac{f'}{f}\varphi'+\beta\frac{f''}{f}\varphi\leq0.
	\end{equation}
	Therefore $\Delta\varphi\leq\beta K\varphi$ distributionally.
\end{lemma}
\begin{proof}
	The differentiability of $f$ and $\varphi$ at $r_0,R$ are straightforward from definitions. Note that $f|_{[r_1,r_0]}$ can be written as
	\begin{equation}\label{eq-countereg:formula_f}
		f(r)=r_0^{-p}\frac{\sqrt{1+c_3^2}}{c_3}\sin\Big[-\sqrt A(r-r_0)+\tan^{-1}(c_3)\Big].
	\end{equation}
	It is therefore clear that $f'(r_1^+)=1$. The fact that $\varphi>0$ follows from $r_0>2/\sqrt{\beta A}$. Then one computes
	\[\Big(f\varphi''+f'\varphi'+\beta f''\varphi\Big)\Big|_{[R,R+c_1R]}\leq R^{1-p}\big[-\frac1{c_1}(1-\beta p)+\frac12(1+\beta)p\big]<0,\]
	\[\Big(f\varphi''+f'\varphi'+\beta f''\varphi\Big)\Big|_{[r_0,R]}=2r^{-p-1}(1-\beta)^2(2\beta-1)<0\]
	To verify (\ref{eq-countereg:supersol}) in $[r_1,r_0]$, note that $-\sqrt A(r-r_0)+\tan^{-1}(c_3)\in[\tan^{-1}(c_3),\pi-\cos^{-1}(c_2)]$, hence $|f'/f|\leq\frac34\sqrt A$ by the definition of constants. Moreover, $\varphi(r)$ has the form $\varphi=a\cosh(\sqrt{\beta A}(r-b)/2)$ for some $a>0,b$. Hence
	\[\Big(\frac{\varphi''}{\varphi}+\frac{f'}{f}\frac{\varphi'}\varphi+\beta\frac{f''}f\Big)\Big|_{[r_1,r_0]}\leq\frac{\beta A}4+\frac34\frac{\sqrt{\beta}A}{2}-\beta A<0.\]
	It remains to verify $\varphi'(r_1^+)<0$. This is equivalent to
	\[\tanh\frac{\sqrt{\beta A}}{2}(r_1-r_0)<-\frac{2}{r_0\sqrt{\beta A}}\]
	The left hand side is smaller than $-\tanh\frac14(\frac\pi2-\cos^{-1}(c_2))\approx -0.16$, while the right hand side is equal to $-\frac 2{c_3p\sqrt\beta}>-\frac 4{c_3}=-0.1$. This proves $\varphi'(r_1^+)<0$.
\end{proof}

\vspace{12pt}

\noindent\textit{{Kai Xu,}}

\vspace{2pt}

\noindent\textit{{Department of Mathematics, Duke University, Durham, NC 27708, USA,}}

\vspace{2pt}

\noindent\textit{Email address: }\url{kai.xu631@duke.edu}.

\end{document}